\documentclass[10pt,reqno]{amsart}
\usepackage[T1]{fontenc}
\usepackage{lmodern}
\usepackage{color}
\usepackage{amsmath,amssymb,mathtools}
\usepackage{mathrsfs}
\usepackage{needspace}
\usepackage{enumerate}
\usepackage[expansion=false]{microtype}
\usepackage[hidelinks]{hyperref}
\numberwithin{equation}{section}
\newtheorem{theorem}{Theorem}[section]
\newtheorem{proposition}[theorem]{Proposition}
\newtheorem{lemma}[theorem]{Lemma}
\newtheorem{corollary}[theorem]{Corollary}
\newtheorem*{problem}{Problem}
\theoremstyle{definition}

\theoremstyle{remark}

\newcommand{\R}{\mathbb R}
\newcommand{\N}{\mathbb N}
\newcommand{\Z}{\mathbb Z}
\newcommand{\E}{\mathbb E}
\newcommand{\D}{\mathcal D}
\newcommand{\one}{\mathbf1}
\newcommand{\PP}{\mathbb P}
\newcommand{\TT}{\mathbb T}
\newcommand{\Tiles}{\mathbb A}
\newcommand{\VV}{\mathbb V}
\newcommand{\dd}{\,\mathrm d}

\newcommand{\ip}[2]{\langle #1,#2\rangle}
\newcommand{\norm}[1]{\lVert #1\rVert}

\newcommand{\Md}{M_{\!d}}
\DeclareMathOperator{\size}{size}
\DeclareMathOperator{\dens}{dens}
\DeclareMathOperator{\Var}{V}

\DeclareMathOperator{\Int}{int}
\newcommand{\cvx}{\mathrm{cvx}}
\newcommand{\Wq}{\mathrm W_q}

\title[Variation Carleson theorem in UMD spaces]{A variation norm Carleson theorem for Walsh--Fourier series in UMD spaces}
\date{}

\author[J. Chang]{Jie Chang}
\address{School of Mathematics and Statistics, Central South University, Changsha 410085,
	People's Republic of China}
\email{changjie@csu.edu.cn}

\author[Y. Jiao]{Yong Jiao}
\address{School of Mathematics and Statistics, Central South University, Changsha 410085,
	People's Republic of China}
\email{jiaoyong@csu.edu.cn}

\author[A. Os\k{e}kowski]{Adam Os\k{e}kowski}
\address{Faculty of Mathematics, Informatics and Mechanics, University of Warsaw,
	Banacha 2, 02-097 Warsaw, Poland}
\email{A.Osekowski@mimuw.edu.pl}

\author[L. Wu]{Lian Wu}
\address{School of Mathematics and Statistics, Central South University, Changsha 410085,
	People's Republic of China}
\email{wulian@csu.edu.cn}

\thanks{This work is supported by the NSFC  (Nos. 12125109, W2411005, 12361131578), Narodowe
Centrum Nauki (Poland) (grant 2023/48/Q/ST1/00013, the Natural Science Foundation of Hunan Province (Nos.  2025ZYJ002, 2024JJ1010) and the Science and Technology Innovation Program of Hunan Province (No. 2024RC1016)}


\begin{document}
\begin{abstract}
We establish a variation norm Carleson theorem for Walsh--Fourier series of functions taking values in an arbitrary complex UMD space $X$. In particular, for every $1<p<\infty$, the Walsh--Fourier series of each $f\in L^p([0,1);X)$ converges almost everywhere to $f$ in the norm of $X$. This answers affirmatively the long-standing open problem asked by Rubio de Francia \cite[Problem~4]{Rubio} in the Walsh-Fourier series setting (the trigonometric analogue will be addressed in a forthcoming paper \cite{CJWZ2026}). A key ingredient of our proof is a coefficient inequality for Walsh wave packets associated with pairwise disjoint tiles, which we derive from martingale cotype by means of a two-point convexity inequality. Combined with time-frequency analysis and martingale variation inequalities, this yields the desired variational estimate without recourse to interpolation with a Hilbert space. We also show that finite Walsh tile type, a geometric property introduced in the study of vector-valued Carleson operators,
is equivalent to the UMD property.
\end{abstract}
\maketitle

\section{Introduction}\label{sec:introduction}
The problem of pointwise convergence of Fourier series has played a central role in the development of harmonic analysis. Let $\mathbb T=\R/\Z$ be the one-dimensional torus, equipped with normalized Lebesgue measure. For a scalar-valued integrable function $f$ on $\mathbb T$, consider its trigonometric Fourier partial sums
\[
 S_Nf(x)=\sum_{|n|\le N}\widehat f(n)e^{2\pi i nx},
 \qquad \widehat f(n)=\int_0^1f(t)e^{-2\pi i nt}\dd t.
\]
A fundamental problem, raised by Lusin, was whether these partial sums converge to $f$ almost everywhere whenever $f\in L^2(\mathbb T)$. This problem remained open for several decades until Carleson \cite{Carleson} gave an affirmative answer in 1966. Hunt \cite{Hunt} subsequently extended Carleson's theorem to the whole range $1<p<\infty$. The celebrated Carleson--Hunt theorem can be formulated in terms of the maximal partial sum operator as follows:
\[
 \left\|\sup_{N\ge0}|S_Nf|\right\|_{L^p(\mathbb T)}
 \lesssim_p\norm{f}_{L^p(\mathbb T)},\qquad 1<p<\infty.
\]
The almost everywhere convergence of Fourier series follows from this estimate and the density of trigonometric polynomials in $L^p(\mathbb T)$. The significance of Carleson's work goes well beyond the convergence problem itself. The analysis of the simultaneous localization in time and frequency, inherent in the proof of his theorem and developed further in the subsequent proofs of Fefferman \cite{Fefferman} and Lacey and Thiele \cite{LaceyThiele}, lies at the heart of modern time-frequency analysis. It has proved fundamental in the study of modulation-invariant operators and, in particular, of multilinear singular integrals. We refer to \cite{GrafakosCFA,ThieleBook} for further background on these developments.

A closely related problem concerns Fourier series with respect to the Walsh system. Shortly after Carleson's work, Billard \cite{Billard} established the almost everywhere convergence of Walsh--Fourier series of square-integrable functions. The corresponding result for $L^p$, $1<p<\infty$, was obtained by Sj\"olin \cite{Sjolin}. The Walsh system, besides being of independent interest, provides a dyadic model of the trigonometric system in which wave packets can be localized exactly in both time and frequency. This feature, together with the close connection between Walsh expansions and dyadic martingales, makes the Walsh model particularly useful in the study of the Carleson operator.

A natural question is how far the above results extend to functions taking values in a Banach space $X$. In this context, the UMD  property of $X$ plays a distinguished role. Indeed, the boundedness of the Hilbert transform on $L^p(\R;X)$, as well as the convergence in $L^p(\mathbb T;X)$ of trigonometric Fourier series for all $f\in L^p(\mathbb T;X)$, is equivalent to the UMD property (see \cite{HNVW}). An analogous characterization in terms of the mean convergence of Walsh--Fourier series was obtained by Wenzel \cite{Wenzel}. However, the almost everywhere convergence of vector-valued Fourier series is a considerably more delicate matter. In 1986, Rubio de Francia \cite{Rubio} established the vector-valued Carleson theorem for UMD Banach lattices and raised the question of whether the result remains true for arbitrary UMD spaces. This question is explicitly recorded as Problem~O.4 on p.~496 of the monograph of Hyt\"onen, van Neerven, Veraar and Weis \cite{HNVW}:
\medskip

\begin{problem}
Does the trigonometric Fourier series of every $f\in L^p(\mathbb T;X)$ converge almost everywhere in the norm of $X$ whenever $1<p<\infty$ and $X$ is a UMD space~?
\end{problem}

\medskip

There has been substantial progress on this problem. Weisz \cite{Weisz} established the almost everywhere convergence of Vilenkin--Fourier series under the additional assumption that the UMD space has an unconditional basis (see also \cite[Introduction]{HLWalsh}). A significant advance was made by Hyt\"onen and Lacey in \cite{HLFourier,HLWalsh}, where the vector-valued Carleson theorem and its Walsh analogue were established for intermediate UMD spaces, that is, spaces of the form $X=[Y,H]_\theta$, where $Y$ is a UMD space, $H$ is a Hilbert space and $0<\theta<1$. In their work, a central role is played by the notion of tile type, a geometric property introduced to replace the orthogonality arguments available in the scalar setting. The question of whether every UMD space admits such a finite tile type, and consequently whether the Walsh--Carleson theorem holds for arbitrary UMD spaces, remained an essential difficulty in this approach.

The primary objective of the present paper is to establish the Walsh analogue of the above problem in full generality. Actually, we will prove a stronger result, concerning the variation of the Walsh--Fourier partial sums. While a maximal inequality controls the magnitude of individual partial sums, a variational inequality controls their successive oscillations and thus provides more information on their convergence. The study of such inequalities for martingales goes back to the work of L\'epingle \cite{Lepingle}, and their connection with the geometry of Banach spaces was developed by Pisier and Xu \cite{PisierXu}. In the context of Fourier analysis, Oberlin, Seeger, Tao, Thiele and Wright \cite{OSTTW} established a variational strengthening of the Carleson--Hunt theorem: the $r$-variation of the trigonometric Fourier partial sums is bounded on $L^p$ for $1<p<\infty$ and $r>\max\{2,p'\}$, where $p'=p/(p-1)$. A weighted version in the Walsh setting was obtained by Do and Lacey \cite{DoLaceyWalsh}. In the vector-valued setting, Hyt\"onen, Lacey and Parissis \cite{HLP} proved a variation norm Carleson theorem for Walsh--Fourier series under a finite tile-type hypothesis; for intermediate UMD spaces they also obtained an improved range of the exponents. More recently, Amenta and Uraltsev \cite{AUVariation} established variational estimates for Fourier integrals with values in intermediate UMD spaces. Our main result dispenses with the extra geometric hypothesis in the Walsh setting.

To state our results, let us fix some notation. Let $(w_n)_{n\ge0}$ be the Walsh system in Paley order
on $[0,1)$. For a Banach space $X$ and
$f\in L^1([0,1);X)$, define $S_0f=0$ and
\begin{equation*}\label{eq:partial}
\begin{aligned}
S_Nf=\sum_{n=0}^{N-1}\widehat f(n)w_n,\qquad \widehat f(n)=\int_0^1 f(x)w_n(x)\dd x
\end{aligned}
\end{equation*}
For $1\le r<\infty$, we define the $r$-variation seminorm
of an $X$-valued sequence $(x_N)_{N\ge0}$  by
\begin{equation*}\label{eq:variation-sequence}
\Var^r(x_N:N\ge0)
=\sup_{K\ge1}\sup_{0\le N_0<\cdots<N_K}
\left(
\sum_{j=1}^K
\norm{x_{N_j}-x_{N_{j-1}}}_X^r
\right)^{1/r}.
\end{equation*}
Our first main result reads as follows; see the next section for any unexplained terminologies.

\begin{theorem}\label{thm:main}
Let $X$ be a complex UMD space. Then there exists
$q_X\in[2,\infty)$ such that, for every $1<p<\infty$
and $r$ satisfying
$\max\{q_X,p'(q_X-1)\}<r<\infty$,
\begin{equation*}
\norm{\Var^r(S_Nf:N\ge0)}_{L^p([0,1))}
\le C_{p,r,X}\norm{f}_{L^p([0,1);X)},
\quad f\in L^p([0,1);X).
\end{equation*}
\end{theorem}

We should mention that Theorem~\ref{thm:main} remains valid
for bounded Vilenkin systems. The proof proceeds almost along
the same lines as in the Walsh case, with appropriate
modifications. To avoid unduly lengthening the paper, we shall not
present the details here. 

Since $S_0f=0$, the maximal partial sum operator is pointwise dominated by $\Var^r(S_Nf:N\ge0)$.
Thus Theorem~\ref{thm:main} immediately leads to following analogue of Carleson--Hunt theorem for
arbitrary UMD spaces.

\begin{corollary}\label{cor:convergence}
Let $X$ be a complex UMD space and $1<p<\infty$.
Then
\begin{equation*}
\left\|\sup_{N\ge0}\norm{S_Nf}_X\right\|_{L^p([0,1))}
\le C_{p,X}\norm{f}_{L^p([0,1);X)},
\qquad f\in L^p([0,1);X).
\end{equation*}
\end{corollary}

Applying Theorem \ref{thm:main}, we may characterize UMD spaces in terms of the almost
everywhere convergence of Walsh--Fourier series, as stated below.

\begin{corollary}\label{cor:walsh-UMD-characterization}
Let $X$ be a complex Banach space. Then $X$ is UMD
if and only if, for some (or equivalently, every)
$1<p<\infty$ and every $f\in L^p([0,1);X)$, $S_Nf(x)\longrightarrow f(x)$
for almost every $x\in[0,1)$.
\end{corollary}

It is worth mentioning that the trigonometric analogues
of Corollaries~\ref{cor:convergence} and
\ref{cor:walsh-UMD-characterization} also hold.
These results and their proofs will be presented in the
forthcoming paper \cite{CJWZ2026}.

Our second main result concerns the geometry of Banach spaces underlying the
vector-valued Walsh--Carleson theorem. It was shown in
\cite[Proposition~4.6]{HLP} that finite Walsh tile type implies the UMD property.
We establish the converse, thus characterizing UMD spaces in terms of Walsh tile type.
Here Walsh tile type is understood in the sense of \cite[Section~4.1]{HLP}; see Subsection~\ref{subsec:tile-type} for the definition.

\begin{theorem}\label{cor:tile-type}
A complex Banach space $X$ is UMD if and only if
it has Walsh tile type $\tau$ for some $2\le\tau<\infty$.
\end{theorem}

Let us briefly explain the main ideas of our proof. A crucial step in extending the classical time-frequency argument to general UMD spaces is to find a suitable substitute for the $L^2$ orthogonality of Walsh wave packets. To this end, given $2\le q<\infty$, we consider the following inequality, referred to as $(\mathrm W_q)$ in Section~\ref{sec:coefficient}:
\[
 \left(\sum_{Q\in\Tiles}\norm{\ip f{w_Q}w_Q}_{L^q(\R_+;X)}^q\right)^{1/q}
 \le A_q\norm f_{L^q(\R_+;X)},
\]
where $\Tiles$ is any finite family of pairwise disjoint tiles in the Walsh phase plane and $w_Q$ is the $L^2$-normalized wave packet associated with $Q$. Here $A_q$ is independent of $f$ and $\Tiles$. An inequality of this nature was established by Amenta and Uraltsev \cite[Proposition~3.8]{AU} for intermediate spaces by complex interpolation. Our approach is different. We first prove that $(\mathrm W_q)$ with constant $1$ is equivalent to the following two-point inequality:
\[
 \left\|\frac{a+b}{2}\right\|_X^q+\left\|\frac{a-b}{2}\right\|_X^q
 \le\frac{\norm a_X^q+\norm b_X^q}{2},\qquad a,b\in X.
\]
The proof relies on compatible subdivisions of dyadic rectangles and the relations between the Walsh wave packets associated with these subdivisions. Using the martingale renorming theorem of Pisier \cite{Pisier,Pisier75}, we then show that every UMD space satisfies $(\mathrm W_q)$ for some finite $q$. This is the principal geometric ingredient of the proof of Theorem~\ref{thm:main}.

To proceed from $(\mathrm W_q)$ to the variational estimate, we combine the time-frequency approach of Hyt\"onen, Lacey and Parissis \cite{HLP} with the reconstruction and coefficient selection arguments of Amenta and Uraltsev \cite{AU}. In particular, the disjoint-packet inequality gives the local packing estimate for the top intervals in the size selection. On individual trees, the UMD property and martingale variation inequalities furnish the required estimates. Combining these estimates by a joint size and density decomposition and interpolation, we obtain Theorem~\ref{thm:main}. As a further consequence of the finite-model estimates, we establish Theorem~\ref{cor:tile-type}. We emphasize that the argument does not require the space to be obtained by complex interpolation with a Hilbert space.

The paper is organized as follows.
Section~\ref{sec:prelim} contains the necessary background
on Banach-valued martingales, the Walsh system and the
dyadic phase plane.
In Section~\ref{sec:coefficient} we establish
$(\mathrm W_q)$ for UMD spaces.
Section~\ref{sec:reduction} reduces the proof of
Theorem~\ref{thm:main} to a uniform $L^p$ bound for
certain linearization operators and also provides the
proof of Corollary~\ref{cor:walsh-UMD-characterization}.
Section~\ref{sec:trees} is devoted to the decomposition
of finite bitile collections into trees, including
the size and density selection and the associated
packing estimates.
In Section~\ref{sec:full} we establish the desired
$L^p$ bound for the linearization operators, thereby
completing the proof of Theorem~\ref{thm:main}.
Finally, Section~\ref{sec:consequences} is devoted to
the proof of Theorem~\ref{cor:tile-type}.

Throughout, all Banach spaces   are complex unless stated otherwise and all vector-valued integrals are understood  in the Bochner sense. We always use $p'=p/(p-1)$ to denote the conjugate exponent of $1<p<\infty$ and $\R_+=[0,\infty)$. The letter $C$ denotes an absolute constant which may change from line to line and $C_\delta$ denotes a constant depending on the parameters $\delta$. We use $\mathbin{\dot\cup}$ and $\mathop{\dot\bigcup}$ to denote disjoint unions.
\medskip

\section{Preliminaries}\label{sec:prelim}\label{sec:phase}
In this section, we collect the conventions and basic facts concerning Banach-valued functions, dyadic analysis and the phase plane that will be used throughout the paper. In particular, we record the relations among Walsh packets, Haar functions and dyadic conditional expectations, together with the order structure of bitiles and trees.

\subsection{Banach-valued functions and the UMD property}\label{subsec:UMD-definition}
For a measure space $(S,\mu)$ and $1\le a<\infty$, the space $L^a(S;X)$ consists of strongly measurable functions $f:S\to X$ for which
\[
 \norm f_{L^a(S;X)}
   =\left(\int_S\norm{f(s)}_X^a\dd\mu(s)\right)^{1/a}<\infty,
\]
with functions identified almost everywhere. For $a=\infty$ we use the essential supremum of the $X$-norm. When the measure space is clear, we write $L^a(X)$ for $L^a(S;X)$. We use the duality convention $\ip{x}{x^*}=x^*(x)$ for complex-linear $x^*\in X^*$.

Let $(\Omega,\mathscr F,\mu)$ be a probability space with a filtration
$\mathscr F_0\subseteq\cdots\subseteq\mathscr F_n$. For $1\le a<\infty$, we write $(x_k)_{k=0}^n$ for an adapted $X$-valued $L^a$ martingale, so that $x_k=\E(x_n\mid\mathscr F_k)$. Its differences $d_k=x_k-x_{k-1}$ satisfy
\[
 d_k\in L^a(\Omega,\mathscr F_k;X),\qquad
 \E(d_k\mid\mathscr F_{k-1})=0\quad(1\le k\le n).
\]
The space $X$ has the UMD property if, for some $1<a<\infty$ and $C_{a,X}<\infty$, these differences satisfy the uniform martingale-transform inequality
\begin{equation}\label{umd def}
 \left\|\sum_{k=1}^n\varepsilon_kd_k\right\|_{L^a(\Omega;X)}
 \le C_{a,X}\left\|\sum_{k=1}^n d_k\right\|_{L^a(\Omega;X)}.
\end{equation}
Here $n\ge1$ and the signs $\varepsilon_k\in\{-1,1\}$ are arbitrary. If the inequality holds at one exponent, it holds at every $1<a<\infty$; see \cite[Definition~8.7 and Theorem~8.12]{Pisier}.

%

\subsection{Dyadic intervals and the Walsh system}\label{subsec:dyadic-system}
We write $\N=\{0,1,2,\ldots\}$ and work on $\R_+=[0,\infty)$ with the dyadic grid
\[
 \D=\bigl\{[m2^{-k},(m+1)2^{-k}):k\in\Z,\ m\in\N\bigr\}.
\]
For $I=[a,b)\in\D$, set
\[
 |I|=b-a,\qquad c(I)=\frac{a+b}{2},\qquad
 I_-=[a,c(I)),\qquad I_+=[c(I),b).
\]
Every $I\in\D$ has a unique dyadic parent
$I^{(1)}\in\D$ of length $2|I|$. If $I=[m2^{-k},(m+1)2^{-k})$, its parent is
\[
 I^{(1)}=[\lfloor m/2\rfloor2^{1-k},(\lfloor m/2\rfloor+1)2^{1-k}).
\]
For frequency intervals we use both half-open conventions. If $\omega=[a,b)$,  write
\[
 [\omega)=[a,b),\qquad (\omega]=(a,b],\qquad \Int\omega=(a,b).
\]

For $n\in\mathbb N$, write $n=\sum_{j\ge0}n_j2^j$, where $n_j\in\{0,1\}$. The Rademacher functions and the Walsh functions in Paley order are given by
\[
 \begin{aligned}
 r_j(t)=(-1)^{\lfloor2^{j+1}t\rfloor},\qquad
 w_n(t)=\prod_{j\ge0}r_j(t)^{n_j}.
 \end{aligned}
\]
Here $j,n\in\N$, $t\in\R_+$, and the product contains only finitely many nontrivial factors. Multiplication of Walsh functions corresponds to addition of their binary digits modulo two; namely, if
\[
 n\oplus m
   :=\sum_{j\ge0}\bigl((n_j+m_j)\bmod2\bigr)2^j,
\]
then $w_nw_m=w_{n\oplus m}$.

For $k\in\Z$, let
$\mathscr D_k=\sigma\{I\in\D:|I|=2^{-k}\}$.
For $f\in L^1(\R_+;X)$, we denote the corresponding conditional expectation
by $\E_k$:
\[
 \E_kf=\E(f\mid\mathscr D_k),\qquad
 \E_kf(x)=\sum_{\substack{I\in\D\\|I|=2^{-k}}}
              \one_I(x)\frac1{|I|}\int_I f(y)\dd y.
\]
The averaging formula gives a representative that is constant on each
interval of length $2^{-k}$, and it also defines $\E_kf$ for
$f\in L^1_{\mathrm{loc}}(\R_+;X)$.
For $1\le p\le\infty$, this operator is a contraction on $L^p(\R_+;X)$.
The tower identity is $\E_j\E_k=\E_{\min\{j,k\}}$ for $j,k\in\Z$. It is well-known that the Walsh functions $(w_n)_{0\le n<2^m}$ form an
orthonormal basis for the space of functions constant on
dyadic intervals of length $2^{-m}$. Consequently,
\[
   \mathbb E_m f
    = \sum_{n=0}^{2^m-1}\widehat f(n)w_n
    = S_{2^m}f,
    \qquad \forall\ f\in L^1([0,1);X).
\]

The differences of successive dyadic conditional expectations can be expressed by Haar expansions.  More precisely, for $I\in\D$, we define the Haar function
\[
 h_I:=|I|^{-1/2}(\one_{I_-}-\one_{I_+}).
\]
Then the function $h_I$ is constant on each child of $I$ and has integral zero.
Therefore, for $j\in\Z$,
\[
 \E_jh_I=
 \begin{cases}
 0,&2^{-j}\ge|I|,\\
 h_I,&2^{-j}<|I|.
 \end{cases}
\]
If $h^{(k)}$ is the sum of the Haar terms of a finite Haar expansion $h$
on intervals of length $2^{-k}$, then
\[
 h^{(k)}=\E_{k+1}h-\E_kh.
\]

\subsection{Tiles and wave packets}
We work in the Walsh phase plane $\R_+\times\R_+$, with time in the first coordinate and frequency in the second. As in \cite[Section~2]{HLWalsh} and \cite[Section~3]{HLP}, we call dyadic rectangles of area one tiles and those of area two bitiles. Their respective collections are
\[
 \begin{aligned}
 \mathcal Q&=\{I\times\omega:I,\omega\in\D,\ |I||\omega|=1\},\\
 \mathcal B&=\{I\times\omega:I,\omega\in\D,\ |I||\omega|=2\}.
 \end{aligned}
\]
For $Q=I_Q\times\omega_Q\in\mathcal Q$, there is a unique $n\in\N$ with
\[
 \omega_Q=|I_Q|^{-1}[n,n+1),\qquad
 w_Q(x)=|I_Q|^{-1/2}\one_{I_Q}(x)w_n(x/|I_Q|).
\]
Thus $\norm{w_Q}_2=1$, $|w_Q|=|I_Q|^{-1/2}\one_{I_Q}$, and, for $f\in L^1_{\mathrm{loc}}(\R_+;X)$,
\[
 \ip f{w_Q}:=\int_{I_Q}f(x)w_Q(x)\dd x\in X.
\]
For $P=I_P\times\omega_P\in\mathcal B$, write
$\omega_P=|I_P|^{-1}[2n,2n+2)$. Its frequency components are
\[
 \begin{aligned}
 \omega_{P_d}&=|I_P|^{-1}[2n,2n+1),& P_d&=I_P\times\omega_{P_d},\\
 \omega_{P_u}&=|I_P|^{-1}[2n+1,2n+2),& P_u&=I_P\times\omega_{P_u}.
 \end{aligned}
\]
In particular, $P=P_d\mathbin{\dot\cup}P_u$ and $P_d,P_u\in\mathcal Q$.

A dyadic rectangle $R=I\times\Omega$ admits a partition into tiles
if and only if $|R|\ge1$.
Indeed, if $\mathscr Q\subset\mathcal Q$ is such a partition, then
\[
 |R|=\sum_{Q\in\mathscr Q}|Q|=\#\mathscr Q\in\{1,2,\ldots\}.
\]
Conversely, if $|R|\ge1$, write $|R|=2^m$ with $m\in\N$.
Dividing $I$ into $2^m$ equal dyadic intervals $J$ gives the tile partition
$\{J\times\Omega\}$.
In particular, if $1\le|R|<2$, every such partition consists of the single
tile $R$.
A dyadic rectangle $P=I\times\omega_P$ of area two can be split into tiles in two ways:
\[
 P=P_d\mathbin{\dot\cup}P_u
  =(I_-\times\omega_P)\mathbin{\dot\cup}(I_+\times\omega_P).
\]
Indeed, for $P=I\times |I|^{-1}[2n,2n+2)$, we have the identities
\[
 w_{2n}(t)=w_n(2t),\qquad w_{2n+1}(t)=r_0(t)w_n(2t).
\]
Therefore,
\begin{equation}\label{eq:bitile-basis-change}
 \begin{pmatrix}w_{P_d}\\ w_{P_u}\end{pmatrix}
 =\frac1{\sqrt2}
   \begin{pmatrix}1&1\\1&-1\end{pmatrix}
   \begin{pmatrix}w_{I_-\times\omega_P}\\w_{I_+\times\omega_P}\end{pmatrix}.
\end{equation}
Thus the time and frequency subdivisions give two orthonormal bases of the same two-dimensional space, as in the recursive construction of \cite[Lecture~5]{ThieleLectures}. For example, the partial-sum representation in \cite[Lecture~6]{ThieleLectures} assigns to $P$ the term
$\ip f{w_{P_d}}w_{P_d}(x)\one_{\{N\in[\omega_{P_u})\}}$.

Throughout, disjointness of tiles refers to the time--frequency plane. In particular, disjoint tiles may have nested time intervals, but their frequency intervals must then be disjoint.

For $I,J\in\D$ with $I\subseteq J$ and $|J|=2^d|I|$, restriction of a Walsh function gives
\begin{equation}\label{eq:restriction}
 \one_I(x)w_n(x/|J|)
   =\epsilon_{I,J,n}\one_I(x)w_{\lfloor n/2^d\rfloor}(x/|I|),
 \qquad \epsilon_{I,J,n}\in\{-1,1\}.
\end{equation}
Indeed, writing $n=\sum_{j\ge0}n_j2^j$, we have
\[
 w_n(x/|J|)
  =\underbrace{\prod_{j<d}r_j(x/|J|)^{n_j}}_{\text{constant on }I}
       \prod_{j\ge d}r_{j-d}(x/|I|)^{n_j}.
\]
The first product equals
$\epsilon_{I,J,n}\in\{-1,1\}$ on $I$ and the second
is $w_{\lfloor n/2^d\rfloor}(x/|I|)$.

The wave packets of disjoint tiles are orthogonal. If $I_Q\cap I_{Q'}=\varnothing$, their packets have disjoint supports. Otherwise assume $I_Q\subseteq I_{Q'}$ and put
\[
 n=|I_Q|\inf\omega_Q,\qquad
 n'=\left\lfloor |I_Q|\inf\omega_{Q'}\right\rfloor.
\]
The interval $|I_Q|^{-1}[n',n'+1)$ contains $\omega_{Q'}$. Disjointness of $\omega_Q$ and $\omega_{Q'}$ implies $n\ne n'$. Therefore, it follows \eqref{eq:restriction} and scalar Walsh orthogonality that
\[
 \begin{aligned}
 \ip{w_Q}{w_{Q'}}
  &=\frac{\pm1}{|I_Q|^{1/2}|I_{Q'}|^{1/2}}
       \int_{I_Q}w_n(x/|I_Q|)w_{n'}(x/|I_Q|)\dd x=0.
 \end{aligned}
\]

\subsection{The phase-plane order and trees}\label{subsec:tile-orders}\label{subsec:tree-geometry}
For two tiles, or two bitiles, we write $P\le T$ precisely when $I_P\subseteq I_T$ and $\omega_T\subseteq\omega_P$. For bitiles, $P\le_uT$ means that $I_P\subseteq I_T$ and $\omega_{T_u}\subseteq\omega_{P_u}$; similarly, $P\le_dT$ means that $I_P\subseteq I_T$ and $\omega_{T_d}\subseteq\omega_{P_d}$. Both component orders imply $P\le T$. Conversely, if $P\le T$, dyadic nesting implies $P\le_uT$ or $P\le_dT$. We write $P<T$ when $P\le T$ and $P\ne T$.

If two dyadic rectangles of the same area intersect, their time intervals are nested, and their lengths force the reverse inclusion of their intersecting frequency intervals. They are therefore comparable: $P\le T$ or $T\le P$. Consequently, if $P,T\in\mathcal Q$ or $P,T\in\mathcal B$, with $P\ne T$ and neither $P\le T$ nor $T\le P$, then $P\cap T=\varnothing$.

A tree consists of a finite family $\TT\subseteq\mathcal B$ and a chosen top $T\in\mathcal B$ such that
\[
 \TT\subseteq\{P\in\mathcal B:P\le T\}.
\]
The top need not be a member of $\TT$. The tree is an up-tree or a down-tree, respectively, when
\[
 \TT\subseteq\{P:P\le_uT\},\qquad
 \TT\subseteq\{P:P\le_dT\}.
\]
For a general tree, put
\[
 \TT_u:=\{P\in\TT:P\le_uT\},\qquad
 \TT_d:=\TT\setminus\TT_u,\qquad
 \TT=\TT_u\mathbin{\dot\cup}\TT_d.
\]
Then $\TT_d\subseteq\{P:P\le_dT\}$. If $T\in\TT$, then $T\in\TT_u$.

A bitile family $\Tiles\subseteq\mathcal B$ is convex precisely when, for every $P,Q\in\Tiles$ and every bitile $R$ with $P\le R\le Q$, one has $R\in\Tiles$.
For a finite bitile family $\PP$, put
\[
 \PP^{\cvx}:=\{R\in\mathcal B:\exists P,Q\in\PP\text{ with }P\le R\le Q\}.
\]

On a tree, a single Walsh modulation converts the two frequency components into a Haar function and a normalized indicator. Fix $Q_T\in\{T_d,T_u\}$ and set
\[
 m_T=|I_T|^{1/2}w_{Q_T},\qquad m_T^2=\one_{I_T}.
\]
For $P\le T$, the interval $\omega_{Q_T}$ lies in one of the two components of $\omega_P$. By \eqref{eq:restriction}, there is a sign $\eta_P$ such that
\begin{equation}\label{eq:two-components}
 (m_Tw_{P_d},m_Tw_{P_u})=
 \begin{cases}
 \eta_P\bigl(|I_P|^{-1/2}\one_{I_P},h_{I_P}\bigr),
     &\omega_{Q_T}\subseteq\omega_{P_d},\\
 \eta_P\bigl(h_{I_P},|I_P|^{-1/2}\one_{I_P}\bigr),
     &\omega_{Q_T}\subseteq\omega_{P_u}.
 \end{cases}
\end{equation}
Indeed, after restriction to $I_P$, the modulating Walsh index is that of the component containing $\omega_{Q_T}$. Exclusive-or with the indices of the two components then gives $0$ and $1$.
For an up-tree $\TT$ with top $T$, we have
$\omega_{T_u}\subseteq\omega_{P_u}$ for every
$P\in\TT$. Thus, choosing $Q_T=T_u$ in
\eqref{eq:two-components}, we obtain
$
m_Tw_{P_d}=\eta_Ph_{I_P},
\, P\in\TT;
$
compare \cite[Lemma~3.2]{HLP}.

\subsection{Walsh tile type}\label{subsec:tile-type}

Following \cite[Section~4.1]{HLP}, a finite family
$(\TT_j)_{j=1}^M$ of up-trees is called \emph{u-good}
if, after a suitable reordering, there exist
corresponding tops $T_j$ such that, with
$\PP=\bigcup_{j=1}^M\TT_j$,
\begin{equation}
 \TT_j=\{P\in\PP:P\le_u T_j,\ P\nleq T_k\text{ for all }1\le k<j\},
 \qquad 1\le j\le M,\label{eq:good-forest}
\end{equation}
and
 \begin{equation} c(\omega_{T_1})\le\cdots\le c(\omega_{T_M}).\label{eq:good-forest-centers}
 \end{equation}
A finite family $(\TT_j)_{j=1}^M$ of down-trees is called
\emph{d-good} if, after a suitable reordering, there
exist corresponding tops $T_j$ satisfying
\eqref{eq:good-forest} with $P\le_d T_j$ in place of
$P\le_u T_j$, and \eqref{eq:good-forest-centers}
with the inequalities reversed.
In either case, the trees are pairwise disjoint as sets of bitiles.
Let $2\le\tau<\infty$. A Banach space $X$ has
\emph{Walsh tile type} $\tau$ if there exists a constant
$C_\tau<\infty$ such that
\begin{equation}\label{eq:tile-type}
\sum_{j=1}^M
\left\|
\sum_{P\in\TT_j}
\ip f{w_{P_d}}w_{P_d}
\right\|_{L^\tau(\R_+;X)}^\tau
\le C_\tau^\tau\norm{f}_{L^\tau(\R_+;X)}^\tau
\end{equation}
for every $f\in L^\tau(\R_+;X)$ and every
u-good family $(\TT_j)_{j=1}^M$ of up-trees.
By \cite[Lemma~4.1 and Remark~4.2]{HLP},
the same condition is obtained by requiring
\eqref{eq:tile-type} for d-good families
of down-trees, with $w_{P_d}$ replaced by
$w_{P_u}$.

\medskip

\section{The disjoint-packet condition for UMD spaces}\label{sec:coefficient}\label{sec:UMD}

Let $2\le q<\infty$ and $0<A_q<\infty$. We say that $X$ satisfies $(\mathrm W_q)$ with constant $A_q$ if
\begin{equation*}\label{eq:Wq-intro}
\left(\sum_{Q\in\Tiles}\bigl\|\ip f{w_Q}w_Q\bigr\|_{L^q(\R_+;X)}^q\right)^{1/q}
\le A_q \norm f_{L^q(\R_+;X)}, \eqno{(W_q)}
\end{equation*}
for every $f\in L^q(\R_+;X)$ and every finite $\Tiles\subseteq\mathcal Q$ whose tiles are pairwise disjoint; that is, $Q\cap Q'=\varnothing$ for all distinct $Q,Q'\in\Tiles$.

Amenta and Uraltsev \cite[Proposition~3.8, (3.10)]{AU}
proved a disjoint-packet estimate for $r$-Hilbertian spaces
by complex interpolation. Their result is formulated in the $3$-Walsh model, but the argument also applies to the
binary Walsh model \cite[p.~4]{AU}. After converting their $L^1$-normalized packets to the
$L^2$ normalization used here, their estimate is precisely
$(\mathrm W_r)$. The following theorem shows that every complex UMD space
$X$ satisfies $(\mathrm W_{q_X})$ for some finite exponent
$q_X\ge2$ depending on $X$, without a Hilbertian
interpolation assumption.

\begin{theorem}\label{prop:UMD-Wq}
Every complex UMD space satisfies $(\mathrm W_{q_X})$ for some finite exponent $q_X\ge2$.
\end{theorem}

For the sake of clarity, we split the proof into three steps.
First, we use rectangular projections and compatible binary refinements
to establish the equivalence between $(\mathrm W_q)$ with constant $1$
and the two-point inequality.
Second, we use a dyadic martingale cotype estimate to obtain an
equivalent norm satisfying a power-type convexity inequality,
and increase the exponent to obtain the two-point inequality
and hence $(\mathrm W_q)$, with explicit bounds on the exponent and
constant. Finally, we use the UMD property to complete the proof.

\medskip

\noindent\textit{Step 1. A two-point criterion.}

\begin{theorem}\label{thm:unit-equivalence}
	Let $X$ be a Banach space and let $2\le q<\infty$.
	Then $X$ satisfies $(\mathrm W_q)$ with constant 1
	if and only if
	\begin{equation}\label{eq:coround}
		\left\|\frac{a+b}{2}\right\|_X^q
		+\left\|\frac{a-b}{2}\right\|_X^q
		\le\frac{\norm a_X^q+\norm b_X^q}{2},
		\qquad a,b\in X.
	\end{equation}
\end{theorem}
Inequality~\eqref{eq:coround} is the averaging form of
\emph{coroundness $q$}; see \cite[Definition~4.1]{ADR}.
For a complex Banach space, coroundness is understood as a property
of its underlying real normed space.

To prove the above theorem, we need some preparation. Let $R=I\times\Omega$ be a dyadic rectangle of area
$N=2^m$ with $m\in\N$. Write
\[
 I=[\ell_I,\ell_I+|I|),\qquad
 \Omega=|I|^{-1}[a,a+N),\qquad a\in N\N,
\]
and partition $I$ into its $N$ equal dyadic subintervals
\[
I_j=\left[\ell_I+\frac{j|I|}{N},
\ell_I+\frac{(j+1)|I|}{N}\right),
\qquad 0\le j<N.
\]
For $h\in L^1(I;X)$, let $\E_{I,N}h$ be the conditional expectation on $I$
with respect to $\sigma\{I_j:0\le j<N\}$, extended by zero outside $I$:
\[
 \E_{I,N}h
   =\sum_{j=0}^{N-1}\one_{I_j}\frac1{|I_j|}\int_{I_j}h.
\]
Set $m_R=\one_Iw_a(\cdot/|I|)$. For
$f\in L^1_{\mathrm{loc}}(\R_+;X)$, define
\begin{equation*}\label{eq:rectangle-proj}
 \Pi_Rf=m_R\E_{I,N}(m_Rf).
\end{equation*}
Since $m_R^2=\one_I$ and $\E_{I,N}^2=\E_{I,N}$, we have
$\Pi_R^2=\Pi_R$.

For $0\le\ell<N$, the condition $a\in N\N$ implies
\[
 a\oplus\ell=a+\ell,\qquad w_{a+\ell}=w_aw_\ell.
\]
Note that the Walsh functions are $1$-periodic and $\ell_I/|I|\in\N$.
Thus 
\[
 \E_{I,N}(m_Rf)
   =\frac1{|I|}\sum_{\ell=0}^{N-1}
      \one_Iw_\ell(\cdot/|I|)
      \int_I m_R(x)f(x)w_\ell(x/|I|)\dd x.
\]
The functions
\[
 e_\ell(x)=|I|^{-1/2}\one_I(x)w_{a+\ell}(x/|I|),
 \qquad 0\le\ell<N,
\]
form an orthonormal basis of the scalar space consisting of $m_R$
times functions constant on $I_j$; see
\cite[Lecture~5, pp.~3--5, and Lecture~6, Section~1.1]{ThieleLectures}.
Multiplying the expansion of $\E_{I,N}(m_Rf)$ by $m_R$, we get
\[
 \Pi_Rf=\sum_{\ell=0}^{N-1}\ip f{e_\ell}e_\ell.
\]
If $R=Q$ is a tile, then $N=1$ and $e_0=w_Q$. Thus, by the conditional Jensen inequality, we have
\begin{equation}\label{eq:tile-proj}
 \begin{aligned}
 &\Pi_Qf=\ip f{w_Q}w_Q\quad\text{when }|Q|=1,\\
 &\norm{\Pi_Rf}_{L^q(\R_+;X)}\le\norm{f\one_I}_{L^q(\R_+;X)}.
 \end{aligned}
\end{equation}

For $N\ge2$, the time and frequency subdivisions of $R=I\times\Omega$ are
\[
 \begin{aligned}
 \mathcal C_{\mathrm t}(R):=\{I_-\times\Omega,I_+\times\Omega\},\qquad
 \mathcal C_{\mathrm f}(R):=\{I\times\Omega_d,I\times\Omega_u\},
 \end{aligned}
\]
where $\Omega_d=[\inf\Omega,c(\Omega))$ and
$\Omega_u=[c(\Omega),\sup\Omega)$. Each family partitions $R$ into
two rectangles of area $N/2$.

For a time subdivision, write $R_\pm=I_\pm\times\Omega$.
The lower frequency index of $R_\pm$ is $a/2$.
As $a$ is even, $w_a(t)=w_{a/2}(2t)$, so that
$m_{R_\pm}=\one_{I_\pm}m_R$.
The time subdivisions of $R_\pm$ are precisely those $I_j$
contained in $I_\pm$. Hence
\begin{equation}\label{eq:time-split}
	\begin{split}
		&\Pi_{R_\pm}f=\one_{I_\pm}\Pi_Rf,\\
		&\norm{\Pi_{R_-}f}_{L^q(\R_+;X)}^q
		+\norm{\Pi_{R_+}f}_{L^q(\R_+;X)}^q
		=\norm{\Pi_Rf}_{L^q(\R_+;X)}^q.
	\end{split}
\end{equation}
For a frequency split, write $R_d=I\times\Omega_d$ and $R_u=I\times\Omega_u$. For $0\le j<N/2$, let $u_j,v_j\in X$ be the values of
$m_R\Pi_Rf=\E_{I,N}(m_Rf)$ on $I_{2j}$ and $I_{2j+1}$, respectively. The carriers of $R_d$ and $R_u$ are $m_R$ and
$m_Rw_{N/2}(\cdot/|I|)$, respectively.
Since $w_{N/2}(\cdot/|I|)$ takes the values $1,-1$
on each pair $I_{2j},I_{2j+1}$, averaging over these pairs, we have
\[
 \begin{aligned}
 m_R\Pi_{R_d}f
   &=\sum_{j=0}^{N/2-1}\frac{u_j+v_j}{2}
                       (\one_{I_{2j}}+\one_{I_{2j+1}}),\\
 m_R\Pi_{R_u}f
   &=\sum_{j=0}^{N/2-1}\frac{u_j-v_j}{2}
                       (\one_{I_{2j}}-\one_{I_{2j+1}}).
 \end{aligned}
\]
Since $\E_{I,N}$ fixes $m_R\Pi_Rf$, the values $u_j,v_j$ do not
change when $f$ is replaced by $\Pi_Rf$. Together with
$\Pi_R^2=\Pi_R$ and \eqref{eq:time-split}, this proves that for
either pair of children $R_1,R_2$,
\begin{equation*}\label{eq:compatibility}
 \Pi_{R_i}\Pi_Rf=\Pi_{R_i}f,\qquad i=1,2.
\end{equation*}
If \eqref{eq:coround} holds, integration over the children gives
\begin{equation}\label{eq:frequency-split}
 \begin{aligned}
 \norm{\Pi_{R_d}f}_q^q+\norm{\Pi_{R_u}f}_q^q
 &=\frac{2|I|}{N}\sum_{j=0}^{N/2-1}
       \left(\left\|\frac{u_j+v_j}{2}\right\|_X^q+
             \left\|\frac{u_j-v_j}{2}\right\|_X^q\right)\\
 &\le\frac{|I|}{N}\sum_{j=0}^{N/2-1}
                  (\norm{u_j}_X^q+\norm{v_j}_X^q)=\norm{\Pi_Rf}_q^q.
 \end{aligned}
\end{equation}

The following lemma shows that any finite family of pairwise disjoint tiles can be included in a partition obtained by successive dyadic subdivisions.

\begin{lemma}\label{lem:refinement}
Let $R_0$ be a dyadic rectangle of area $2^m$, with $m\in\N$, and let $\Tiles$ be
a finite pairwise disjoint family of tiles contained in $R_0$.
There are partitions $(\mathcal R_j)_{j=0}^m$ of $R_0$ such that
\begin{equation}\label{R_0-partition}
 \mathcal R_0=\{R_0\},\qquad
 \mathcal R_{j+1}=\bigcup_{R\in\mathcal R_j}\mathcal C(R),
 \qquad \Tiles\subseteq\mathcal R_m,
\end{equation}
where the recursion is for $0\le j<m$, $\mathcal C(R)$ is either
$\mathcal C_{\mathrm t}(R)$ or $\mathcal C_{\mathrm f}(R)$, and every
$R\in\mathcal R_j$ has
area $2^{m-j}$.
\end{lemma}

\begin{proof}
For a rectangle $R=I\times\Omega$ of area greater than one, set
\[
 \Tiles_R=\{Q\in\Tiles:Q\subseteq R\}.
\]
If $I_Q\ne I$ for every $Q\in\Tiles_R$, then each $I_Q$ is contained
in one of the dyadic halves of $I$, and we choose a time split.
Otherwise, fix $Q\in\Tiles_R$ with $I_Q=I$. Its frequency interval
is a proper subinterval of $\Omega$, since $|Q|=1<|R|$.
If another tile $Q'\in\Tiles_R$ had $\omega_{Q'}=\Omega$, then
\[
 Q\cap Q'=I_{Q'}\times\omega_Q\ne\varnothing,
\]
contrary to disjointness. Thus every frequency interval in $\Tiles_R$
is contained in one of the dyadic halves of $\Omega$, and we choose
a frequency split. When $\Tiles_R$ is empty, the time split applies.
In either case, every tile in $\Tiles_R$ is contained in exactly one
child. Repeating this construction for $m$ generations gives rectangles
of area one. Each tile in $\Tiles$ is contained in one of these
rectangles and hence equals it. The case $m=0$ is immediate.
\end{proof}
\begin{proof}[Proof of Theorem~\ref{thm:unit-equivalence}]
Suppose first that $(\mathrm W_q)$ holds with constant $1$.
	For $a,b\in X$, take
	\[
	f=a\one_{[0,1/2)}+b\one_{[1/2,1)},\qquad
	Q_0=[0,1)\times[0,1),\qquad
	Q_1=[0,1)\times[1,2).
	\]
	Then
	\[
		\ip f{w_{Q_0}}=\frac{a+b}{2},\qquad
		\ip f{w_{Q_1}}=\frac{a-b}{2}.\]
Moreover,	
		\[\norm{f}_{L^q(\R_+;X)}^q
		=\frac{\norm a_X^q+\norm b_X^q}{2}.
	\]
Since both packets have modulus one on $[0,1)$,  applying ($W_q$), we get the desired estimate \eqref{eq:coround}.

Conversely, suppose \eqref{eq:coround} holds. Fix
$f\in L^q(\R_+;X)$ and a finite pairwise disjoint family
$\Tiles\subseteq\mathcal Q$.
Choose powers of two $L,M$ so that
	\[
	\bigcup_{Q\in\Tiles}Q
	\subseteq R_0=[0,L)\times[0,M),\qquad LM=2^m\ge1.
	\]
Applying Lemma~\ref{lem:refinement}, there are partitions
	$(\mathcal R_j)_{j=0}^m$ of $R_0$ such that \eqref{R_0-partition} holds. Then it follows from \eqref{eq:time-split} and \eqref{eq:frequency-split} that
	\[
	\sum_{R'\in\mathcal R_{j+1}}
	\norm{\Pi_{R'}f}_{L^q(\R_+;X)}^q
	\le\sum_{R\in\mathcal R_j}
	\norm{\Pi_Rf}_{L^q(\R_+;X)}^q,
	\qquad 0\le j<m.
	\]
Since $\Tiles\subseteq\mathcal R_m$, it follows from iteration and \eqref{eq:tile-proj} that
	\[
	\begin{aligned}
		\sum_{Q\in\Tiles}
		\bigl\|\ip f{w_Q}w_Q\bigr\|_{L^q(\R_+;X)}^q
		&\le\sum_{Q\in\mathcal R_m}
		\norm{\Pi_Qf}_{L^q(\R_+;X)}^q\\
		&\le\norm{\Pi_{R_0}f}_{L^q(\R_+;X)}^q\le\norm{f}_{L^q(\R_+;X)}^q,
	\end{aligned}
	\]
which is exactly $(\mathrm W_q)$ with constant $1$.
\end{proof}
\medskip

\noindent\textit{Step 2. From dyadic martingale cotype to $(\mathrm W_q)$.}
We start by recalling an equivalent norm theorem from \cite[Theorem~4.6]{Pisier}. The origin of this result goes back to \cite{Pisier75}. Let $X$ be a complex Banach space. For $N\ge0$, equip
$\Omega_N=\{-1,1\}^N$ with the uniform probability measure
and the coordinate filtration
$\mathscr F_n=\sigma(\omega_1,\ldots,\omega_n)$.
Here $\Omega_0$ is a singleton and $\mathscr F_0$ is trivial.
By a finite dyadic martingale we mean a sequence
\[
x_n=\E(x_N\mid\mathscr F_n),\qquad 0\le n\le N.
\]
Its initial value is constant but need not vanish.
If, for some $2\le s<\infty$ and $D\ge1$,
every such martingale satisfies
\begin{equation}\label{eq:mart-s}
	\sum_{n=1}^N\E\norm{x_n-x_{n-1}}_X^s
	\le D^s\E\norm{x_N}_X^s,
\end{equation}
then we say that $X$ satisfies \textit{dyadic martingale cotype} $s$ with constant $D$.
By choose $q=s$, $\alpha=2^{-1/s}$ and $C=2^{1/s}D$ in \cite[Theorem~4.6]{Pisier}, we have

\begin{lemma}\label{lem:renorm}
Suppose that $X$ satisfies dyadic martingale cotype $s$ with constant $D$.
Then $X$ admits an equivalent complex norm $|\cdot|$ such that
\begin{equation}\label{eq:norm-equivalence}
 2^{-1/s}\norm x_X\le |x|\le\norm x_X
\end{equation}
and
\begin{equation}\label{eq:power-convex}
 |x|^s+\delta|y|^s\le\frac{|x+y|^s+|x-y|^s}{2}
 \qquad \forall\ x,y\in X,
\end{equation}
with $\delta=(2D^s)^{-1}$.
\end{lemma}


\medskip
We now show that increasing the exponent allows the coefficient
$\delta$ in \eqref{eq:power-convex} to be replaced by $1$.
For the qualitative connection between power-type convexity
and finite coroundness, see \cite[Remark~5.4]{ADR}.
The following argument gives an explicit exponent.

\Needspace{8\baselineskip}
\begin{lemma}\label{lem:power}
Suppose a norm $|\cdot|$ satisfies \eqref{eq:power-convex} for some $s\ge2$ and $0<\delta\le1$.
If $k$ is a positive integer with $k\delta\ge1$, then $|\cdot|$ satisfies \eqref{eq:coround} with exponent $q=ks$.
\end{lemma}

\begin{proof}
It suffices to prove
\begin{equation}\label{eq:sufficient}
 \frac{|x+y|^{ks}+|x-y|^{ks}}2\ge|x|^{ks}+|y|^{ks}.
\end{equation}
Indeed, applying the above estimate to $(a+b)/2$ and $(a-b)/2$ gives \eqref{eq:coround} with $q=ks$. Note that both sides of \eqref{eq:sufficient} are symmetric in $x,y$, since $|y-x|=|x-y|$, so we may
assume $|x|\ge|y|$. If $x=0$, then $y=0$ and the assertion is immediate. Otherwise, it follows from Jensen's inequality, \eqref{eq:power-convex} and the binomial inequality that
\[
 \begin{aligned}
 \frac{|x+y|^{ks}+|x-y|^{ks}}2
 &\ge\left(\frac{|x+y|^s+|x-y|^s}{2}\right)^k\\
 &\ge(|x|^s+\delta|y|^s)^k\\
 &\ge|x|^{ks}+k\delta|x|^{s(k-1)}|y|^s\\
 &\ge|x|^{ks}+|y|^{ks}.
 \end{aligned}
\]
where we use $|y|\le|x|$ and $k\delta\ge1$ in the last inequality. This verifies the desired assertion \eqref{eq:sufficient} and the proof is complete.
\end{proof}

\begin{corollary}\label{cor:martingale-Wq}
Suppose that $X$ satisfies dyadic martingale cotype $s$ with constant $D$. Then $X$ satisfies $(\mathrm W_q)$ with
\( q=s\lceil2D^s\rceil\) and \(A_q\le2^{1/s}\).

\end{corollary}
\begin{proof}

Take the norm $|\cdot|$ supplied by Lemma~\ref{lem:renorm} and set $\delta=(2D^s)^{-1},$
	$k=\lceil2D^s\rceil$, $q=ks$. Since $k\delta\ge1$,
	Lemma~\ref{lem:power} and Theorem~\ref{thm:unit-equivalence}
	give $(\mathrm W_q)$ with constant $1$ in the norm $|\cdot|$. For every finite family $\Tiles$ of pairwise disjoint tiles and every $f\in L^q(\R_+;X)$, we deduce from Theorem~\ref{thm:unit-equivalence} and \eqref{eq:norm-equivalence} that
\[
 \begin{aligned}
 \sum_{Q\in\Tiles}\int_{\R_+}
       \norm{\ip f{w_Q}w_Q(t)}_X^q\dd t
 &\le 2^{q/s}\sum_{Q\in\Tiles}\int_{\R_+}
       |\ip f{w_Q}w_Q(t)|^q\dd t\\
 &\le 2^{q/s}\int_{\R_+}|f(t)|^q\dd t\le 2^{q/s}\norm f_{L^q(\R_+;X)}^q.
 \end{aligned}
\]
Taking $q$th roots, the desired assertion follows.
\end{proof}

\medskip
\noindent\textit{Step 3. Completion of the proof.}
Since $X$ is UMD, it has finite Rademacher cotype; see
\cite[Theorem~4.40(ii) and Propositions~8.8 and~8.27]{Pisier}.
Thus there exist $2\le s<\infty$ and $C_s\ge1$ such that
\begin{equation}\label{eq:rad-cotype}
    \sum_{k=1}^n\norm{x_k}_X^s
    \le C_s^s\,\E_\varepsilon
    \left\|\sum_{k=1}^n\varepsilon_k x_k\right\|_X^s
\end{equation}
for every finite sequence $(x_k)_{k=1}^n$ in $X$, where
$\E_\varepsilon$ denotes expectation over independent
Rademacher signs. It is well known that, in a UMD space, Rademacher cotype $s$ implies martingale
cotype $s$; see \cite[Proposition~4.37 and its proof]{Pisier}
and \cite[Section~4.2]{HLP}. Below, we recall the short randomization argument to make the
resulting martingale cotype constant explicit. Let $\beta_{s,X}\ge\max\{1,C_{s,X}\}$ be a constant in the $L^s(X)$ UMD
inequality for martingale transforms by deterministic signs;
such a constant exists by \cite[Theorem~8.12]{Pisier}.
Applying \eqref{eq:rad-cotype} pointwise and integrating,
then applying the UMD inequality for each fixed choice of signs,
we obtain
\begin{align*}
    \sum_{k=1}^N\E\norm{d_k}_X^s
    &\le C_s^s\,\E_\varepsilon
    \left\|\sum_{k=1}^N\varepsilon_k d_k\right\|_{L^s(X)}^s
    \notag\\
    &\le (C_s\beta_{s,X})^s
    \norm{x_N-x_0}_{L^s(X)}^s
    \notag\le (2C_s\beta_{s,X})^s
    \norm{x_N}_{L^s(X)}^s.
    \label{eq:UMD-mart}
\end{align*}
In particular, $X$ satisfies dyadic martingale cotype $s$
with constant $D=2C_s\beta_{s,X}$.
Then it follows from Corollary~\ref{cor:martingale-Wq} that $X$  satisfies
	$(\mathrm W_q)$ with
\begin{equation*}\label{eq:explicit-q}
    q_X=s\left\lceil 2D^s\right\rceil,
    \qquad A_{q_X}\le 2^{1/s}.
\end{equation*}
This completes the proof of Theorem \ref{prop:UMD-Wq}.

\medskip

\section{Reduction to  linearization operators}\label{sec:reduction}
This section aims to reduce Theorem~\ref{thm:main} to an $L^p$ estimate for the  linearization operators $L_{\PP}$ defined below. To this end, we first estimate the partial-sum variation in terms of two finite packet models defined on a common family. Reflection then reduces these estimates to the down-packet model, whose value is attained by a measurable  linearization.

\subsection{Reduction of partial-sum variation to finite packet models}\label{subsec:partial-models}
We use finite versions of the packet models introduced in \cite[Section~5]{HLP}. For $K\ge1$, let
\[
 \mathcal N_K=\{(N_0,\ldots,N_K)\in[0,\infty)^{K+1}:
                                      N_0<\cdots<N_K\}.
\]
Fix $1\le r<\infty$, a finite $\PP\subseteq\mathcal B$, and signs
$(\epsilon_P)_{P\in\PP}\in\{-1,1\}^{\PP}$. Define the down-packet model by
\begin{equation}\label{eq:model}
 \begin{split}
 C_{r,\PP}f(x)&=\sup_{K\ge1}\ \sup_{(N_0,\ldots,N_K)\in\mathcal N_K}\\
 &\qquad\Bigg(\sum_{j=1}^K\Bigg\|
 \sum_{P\in\PP}\epsilon_P\ip f{w_{P_d}}w_{P_d}(x)
 \one_{\{\substack{N_j\in[\omega_{P_u})\\
                    N_{j-1}\notin\Int\omega_P}\}}
 \Bigg\|_X^r\Bigg)^{1/r},
 \end{split}
\end{equation}
and the up-packet model by
\begin{equation}\label{eq:symmetric-condition}
 \begin{split}
 \widetilde C_{r,\PP}f(x)
 &=\sup_{K\ge1}\ \sup_{(N_0,\ldots,N_K)\in\mathcal N_K}\\
 &\quad\left(\sum_{j=1}^K\left\|
   \sum_{P\in\PP}\epsilon_P\ip f{w_{P_u}}w_{P_u}(x)
   \one_{\{\substack{N_{j-1}\in(\omega_{P_d}]\\
                      N_j\notin\Int\omega_P}\}}
                          \right\|_X^r\right)^{1/r}
 \end{split}
\end{equation}
for all
$f\in L^1_{\mathrm{loc}}(\R_+;X)$. Both models vanish when $\PP=\varnothing$. For partial sums the endpoints are natural numbers and the packet signs are one. The following proposition is a finite version of the reduction in \cite[Proof of Theorem~1.3]{HLP}, with a common bitile family independent of the endpoint sequence.
%

\begin{proposition}\label{prop:partial-reduction}
For $M\in\N$, let
\[
 \PP_M=\bigl\{P\in\mathcal B:I_P\subseteq[0,1),\
        2^{-M}\le|I_P|\le1,\ \omega_P\subseteq[0,2^{M+1})\bigr\}.
\]
With packet signs equal to one, for every $f\in L^1([0,1);X)$ and $1\le r<\infty$,
\begin{equation}\label{eq:partial-to-model}
 \Var^r(S_Nf:0\le N\le2^M)
 \le C_{r,\PP_M}(f\one_{[0,1)})+\widetilde C_{r,\PP_M}(f\one_{[0,1)}).
\end{equation}
The finite family $\PP_M$ is independent of the endpoint sequence.
\end{proposition}

\begin{proof}
We use the maximal dyadic frequency-block decomposition
in \cite[Proof of Theorem~1.3]{HLP}. For integers
$0\le a<b\le2^M$, let $\mathcal J(a,b)$ be the
family of maximal dyadic intervals contained in
$[a,b)$. These intervals partition $[a,b)$, and
the Walsh projection identity gives
\[
 S_bf-S_af
 =\sum_{J\in\mathcal J(a,b)}
   \sum_{\substack{I_Q\subseteq[0,1)\\
                    \omega_Q=J}}
      \ip f{w_Q}w_Q.
\]
As shown in \cite{HLP}, a block
$J=\omega_{P_d}$ is maximal precisely when
$b\in[\omega_{P_u})$ and
$a\notin\Int\omega_P$; similarly, a block
$J=\omega_{P_u}$ is maximal precisely when
$a\in(\omega_{P_d}]$ and
$b\notin\Int\omega_P$.
These conditions remain valid for $a=0$.
Moreover, every $J\in\mathcal J(a,b)$ satisfies
\[
 1\le |J|\le2^M,\qquad
 J^{(1)}\subseteq[0,2^{M+1}).
\]
Thus all the corresponding bitiles belong to
$\PP_M$. Applying the two model decompositions
to each increment $S_{N_j}f-S_{N_{j-1}}f$
and using the triangle inequality in
$\ell^r(X)$, we obtain the desired estimate \eqref{eq:partial-to-model}.
\end{proof}

%

\subsection{Reduction of the symmetric model to the down-packet model}\label{subsec:modulation}

We show in this subsection that an \(L^p\) estimate for the down-packet model transfers to the up-packet model with the same constant. The proof adapts the frequency-reflection argument of \cite[Lemma~4.1]{HLP}.

\begin{proposition}\label{prop:symmetric-reduction}
Fix $1<p<\infty$ and $1\le r<\infty$. If
$\norm{C_{r,\PP}f}_{L^p}\le C\norm f_{L^p(X)}$
uniformly over all finite bitile families $\mathbb P\subset \mathcal B$ and all choices of signs $(\epsilon_P)_{P\in\mathbb P}\in\{-1,1\}^{\mathbb P}$, then
$\norm{\widetilde C_{r,\PP}f}_{L^p}\le C\norm f_{L^p(X)}$
with the same  constant.
\end{proposition}

\begin{proof}
The case $\PP=\varnothing$ is immediate.
Otherwise choose a positive dyadic number $L$
and $H\in\N$ such that, with $B=2^H/L$,
\[
 I_P\subseteq[0,L),\qquad
 \omega_P\subseteq[0,B),\qquad
 |I_P|=L2^{-d_P},\quad 0\le d_P\le H
 \quad(P\in\PP).
\]
Set
\[
 m(x)=\one_{[0,L)}(x)w_{2^H-1}(x/L).
\]
For every tile or bitile
$R=I_R\times[a,b)$ contained in
$[0,L)\times[0,B)$, define
\[
 R^*=I_R\times[B-b,B-a),\qquad
 \PP^*=\{P^*:P\in\PP\},\qquad
 \epsilon_{P^*}=\epsilon_P.
\]
Frequency reflection exchanges the upper and lower
frequency components. The modulation identity
in \cite[Lemma~4.1]{HLP}, together with
\eqref{eq:restriction}, gives
\[
 mw_Q=\eta_Qw_{Q^*},
 \qquad \eta_Q\in\{-1,1\},
\]
and hence
\[
 \ip{mf}{w_{Q^*}}w_{Q^*}
 =m\ip f{w_Q}w_Q.
\]

The suprema defining both models may be
restricted to endpoints in $[0,B]$:
endpoints larger than $B$ can be replaced
by $B$, and the resulting zero blocks
can be deleted. For
$0\le N_0<\cdots<N_K\le B$, put
\[
 \widetilde N_i=B-N_{K-i},
 \qquad 0\le i\le K.
\]
Then, for $i=K-j+1$,
\[
 \begin{aligned}
 \widetilde N_i\in[\omega_{(P^*)_u})
 &\iff N_{j-1}\in(\omega_{P_d}],\\
 \widetilde N_{i-1}\notin\Int\omega_{P^*}
 &\iff N_j\notin\Int\omega_P.
 \end{aligned}
\]
Therefore, reversing the order of the blocks and using
the modulation identity gives
\[
 \widetilde C_{r,\PP}f=C_{r,\PP^*}(mf).
\]
Applying the assumed estimate to $\mathbb P^*$ and $mf$, we obtain
\[
 \norm{\widetilde C_{r,\PP}f}_{L^p}
 =\norm{C_{r,\PP^*}(mf)}_{L^p}
 \le C\norm{mf}_{L^p(X)}
 \le C\norm f_{L^p(X)},
\]
which proves the desired assertion.
\end{proof}

\subsection{Linearization of the down-packet model}
\label{subsec:linearization}
Fix $1<r<\infty$,
a finite bitile family $\PP\subseteq\mathcal B$,
and signs $\epsilon_P\in\{-1,1\}$ for $P\in\PP$.
Let $K\ge1$ be an integer, let
$
 \mathbf N=(N_0,\ldots,N_K):
 \R_+\longrightarrow\mathcal N_K
$
be a measurable map and let
$a_j^*:\R_+\to X^*$, $1\le j\le K$,
be strongly measurable functions satisfying
\begin{equation}\label{eq:directions}
 \sum_{j=1}^K\norm{a_j^*(x)}_{X^*}^{r'}
 \le1, \qquad x\in\R_+,
\end{equation}
where $r'=\frac r{r-1}$. For $P\in\PP$ and $1\le j\le K$, set
\[
 A(P,j)=\bigl\{x\in\R_+:
 N_j(x)\in[\omega_{P_u}),\
 N_{j-1}(x)\notin\Int\omega_P\bigr\}.
\]
For $f\in L^1_{\mathrm{loc}}(\R_+;X)$, define the scalar linearization operator $L_{\PP}$ by
\begin{equation}\label{eq:linearization}
 L_{\PP}f(x)
 =\sum_{j=1}^K\sum_{P\in\PP}
 \epsilon_Pw_{P_d}(x)
 \ip{\ip f{w_{P_d}}}{a_j^*(x)}
 \one_{A(P,j)}(x),
 \qquad x\in\R_+.
\end{equation}
Then
the operator $L_{\PP}$ is complex-linear.
Moreover, by H\"older's inequality,
\[
 |L_{\PP}f(x)|\le C_{r,\PP}f(x),
 \qquad x\in\R_+.
\]
The next proposition shows that equality
can be attained by suitable measurable choices
of $\mathbf N$ and $a_j^*$.

\begin{proposition}\label{prop:exact-linearization}
Let $1<p<\infty$ and $1<r<\infty$.
For every finite bitile family
$\PP\subseteq\mathcal B$,
$(\epsilon_P)_{P\in\PP}\in\{-1,1\}^{\PP}$
and $f\in L^p(\R_+;X)$,
there exist an integer $K\ge1$,
a measurable map
$\mathbf N:\R_+\to\mathcal N_K$
and strongly measurable functions
$a_j^*:\R_+\to X^*$, $1\le j\le K$,
satisfying \eqref{eq:directions}, such that
the corresponding linearization $L_f$ satisfies
\begin{equation}\label{eq:exact-linearization}
 L_f f=C_{r,\PP}f.
\end{equation}
\end{proposition}

\begin{proof}
If $\PP=\varnothing$, take $K=1$,
$\mathbf N(x)=(0,1)$ and $a_1^*(x)=0$
for all $x\in\R_+$. Then
$L_f f=C_{r,\PP}f=0$.
Suppose $\PP\ne\varnothing$.
For any $\mathbf N=(N_0,\ldots,N_K)\in\mathcal N_K$,
a bitile $P\in\PP$ contributes to at most one block.
Indeed, if $j<j'$ and $P$ contributes to both blocks,
then
\[
 \inf\omega_P<c(\omega_P)
 \le N_j\le N_{j'-1}<N_{j'}<\sup\omega_P,
\]
contrary to $N_{j'-1}\notin\Int\omega_P$.
For $1\le j\le K$, put
\[
 B_j(\mathbf N)
 =\{P\in\PP:N_j\in[\omega_{P_u}),\
 N_{j-1}\notin\Int\omega_P\}.
\]
Thus, the nonempty sets $B_j(\mathbf N)$ form
an ordered list of pairwise disjoint subsets of $\PP$.
Since $\PP$ is finite, only finitely many such lists
can occur. Choose a representative
$\mathbf N^\alpha\in\mathcal N_{K_\alpha}$
for each realizable list, $1\le\alpha\le A$,
retaining any empty blocks in the representative.
Set
\[
 B=\max_{P\in\PP}\sup\omega_P,\qquad
 K_*=\max_{1\le\alpha\le A}K_\alpha.
\]
Extend each representative to length $K_*$ by
\[
 N_j^\alpha
 =\max\{N_{K_\alpha}^\alpha,B\}+j-K_\alpha,
 \qquad K_\alpha<j\le K_*.
\]
The additional blocks select no bitiles, and each
extended sequence remains strictly increasing.

For $1\le\alpha\le A$ and $1\le j\le K_*$, define
\[
 z_{\alpha,j}(x)
 =\sum_{P\in B_j(\mathbf N^\alpha)}
 \epsilon_P\ip f{w_{P_d}}w_{P_d}(x),
\]
where $B_j(\mathbf N^\alpha)=\varnothing$
for $j>K_\alpha$. Since empty blocks contribute
zero, the choice of representatives gives
\[
 C_{r,\PP}f(x)
 =\max_{1\le\alpha\le A}
 \left(\sum_{j=1}^{K_*}
 \norm{z_{\alpha,j}(x)}_X^r\right)^{1/r}.
\]
Choose a dyadic interval $I$ containing all $I_P$,
$P\in\PP$, and partition $I$ into finitely many
dyadic intervals on which every $z_{\alpha,j}$
is constant. For each interval $J$ of this
partition, choose $x_J\in J$ and an index
$\alpha_J$ attaining the above maximum.
For every nonzero $z_{\alpha_J,j}(x_J)$, applying
the Hahn--Banach theorem, there exists
$\phi_{J,j}\in X^*$ such that
\[
 \norm{\phi_{J,j}}_{X^*}=1,\qquad
 \ip{z_{\alpha_J,j}(x_J)}{\phi_{J,j}}
 =\norm{z_{\alpha_J,j}(x_J)}_X.
\]

On each interval $J$ with
$C_{r,\PP}f(x_J)>0$, define
$\mathbf N(x)=\mathbf N^{\alpha_J}$ and
\[
 a_j^*(x)
 =\frac{\norm{z_{\alpha_J,j}(x_J)}_X^{r-1}}
 {C_{r,\PP}f(x_J)^{r-1}}\phi_{J,j},
 \qquad x\in J,\quad 1\le j\le K_*,
\]
with $a_j^*(x)=0$ when
$z_{\alpha_J,j}(x_J)=0$.
On intervals where $C_{r,\PP}f=0$ and outside
$I$, take $\mathbf N(x)=\mathbf N^1$
and $a_j^*(x)=0$.  Then
for every $x\in J$ with $C_{r,\PP}f(x_J)>0$,
$$\sum_{j=1}^{K_*}\norm{a_j^*(x)}_{X^*}^{r'}
 =C_{r,\PP}f(x_J)^{-r}
   \sum_{j=1}^{K_*}\norm{z_{\alpha_J,j}(x_J)}_X^r
 =1,$$
$$
 L_f f(x)
 =\sum_{j=1}^{K_*}
   \ip{z_{\alpha_J,j}(x_J)}{a_j^*(x)}=C_{r,\PP}f(x_J)^{1-r}
   \sum_{j=1}^{K_*}\norm{z_{\alpha_J,j}(x_J)}_X^r
 =C_{r,\PP}f(x).$$
On intervals where \(C_{r,\PP}f=0\) and outside \(I\), both \(L_f f\) and \(C_{r,\PP}f\) vanish.
The map $\mathbf N$ and the functions $a_j^*$
are constant on each interval $J$ of the partition and on $\R_+\setminus I$. Hence $\mathbf N$ is measurable and each $a_j^*$ is strongly measurable.
Thus \eqref{eq:directions} holds and
$L_f f=C_{r,\PP}f$, which gives the desired assertion.
\end{proof}

\subsection{Reduction to the linearization operators}

The following theorem establishes the $L^p$ boundedness
of the linearization operators $L_{\PP}$. Together with Theorem \ref{prop:UMD-Wq},
Propositions~\ref{prop:partial-reduction},
\ref{prop:symmetric-reduction} and
\ref{prop:exact-linearization}, it implies
Theorem~\ref{thm:main}, as shown below.
Its proof is deferred to Section~\ref{sec:full}.

\begin{theorem}\label{thm:fixed-linearization}
Let $X$ be a complex UMD space which satisfies $(\mathrm W_q)$ for some
$2\le q<\infty$. Suppose that $
 1<p<\infty$ and  $\max\{q,p'(q-1)\}<r<\infty.$
For every finite bitile family $\PP\subseteq\mathcal B$, $(\epsilon_P)_{P\in\PP}\in\{-1,1\}^{\PP}$, $K\ge1$, a measurable map $\mathbf N:\R_+\to\mathcal N_K$, and  strongly measurable directions $a_j^*:\R_+\to X^*$ satisfying
\eqref{eq:directions}, the operator $L_{\PP}$  satisfies
\[
 \norm{L_{\PP}h}_{L^p(\R_+)}
 \le C_{p,r,q,A_q,X}\norm h_{L^p(\R_+;X)},
 \qquad\forall\ h\in L^p(\R_+;X).
\]
\end{theorem}

Combining Theorem \ref{thm:fixed-linearization} with Propositions \ref{prop:symmetric-reduction}, \ref{prop:exact-linearization}, we have the following $L^p$ bound for finite-packet models.

\begin{theorem}\label{thm:model}
Under the hypotheses  of
Theorem~\ref{thm:fixed-linearization}, for every
$f\in L^p(\R_+;X)$,
\begin{equation}\label{eq:model-estimate}
 \norm{C_{r,\PP}f}_{L^p(\R_+)}+
 \norm{\widetilde C_{r,\PP}f}_{L^p(\R_+)}
       \le C_{p,r,q,A_q,X}\norm f_{L^p(\R_+;X)}.
\end{equation}
\end{theorem}

\medskip

We  conclude this section the proof of Theorem~\ref{thm:main} and Corollary \ref{cor:walsh-UMD-characterization}.

\begin{proof}[Proof of Theorem~\ref{thm:main}]
Since $X$ is  UMD, by Theorem~\ref{prop:UMD-Wq}, there exist
$q_X\ge2$ and $A_{q_X}>0$
such that $X$ satisfies $(\mathrm W_{q_X})$. Fix $1<p<\infty$ and
$r<\infty$ satisfying  $r>\max\{q_X,p'(q_X-1)\}$. Then Theorem~\ref{thm:model}   applies
with $q=q_X$.  Let $f\in L^p([0,1);X)$ and extend it by zero to $\R_+$. For each $M\in\N$,  by Proposition \ref{prop:partial-reduction} and \eqref{eq:model-estimate}, we obtain
\[
\begin{aligned}
 &\norm{\Var^r(S_Nf:0\le N\le2^M)}_{L^p([0,1))}\\
 &\quad\le
 \norm{C_{r,\PP_M}(f\one_{[0,1)})}_{L^p(\R_+)}
 +\norm{\widetilde C_{r,\PP_M}(f\one_{[0,1)})}_{L^p(\R_+)}\\
 &\quad\le C_{p,r,q_X,X}
             \norm f_{L^p([0,1);X)}.
\end{aligned}
\]
As $M\to\infty$, the finite variations on the left
increase pointwise to $\Var^r(S_Nf:N\ge0)$. We conclude the desired result by the monotone convergence theorem.
\end{proof}

\begin{proof}[Proof of
Corollary~\ref{cor:walsh-UMD-characterization}]
We first show the sufficiency.
Suppose that $X$ is UMD and $1<p<\infty$.
For every $f\in L^p([0,1);X)$, by Theorem~\ref{thm:main}, for sufficiently large $r$
we have
$
 \Var^r(S_Nf:N\ge0)<\infty$
 almost everywhere.
Hence $(S_Nf)_{N\ge0}$ converges almost everywhere
in $X$. Since $S_{2^m}f=\mathbb E_mf\to f$ almost everywhere
by the convergence theorem for dyadic conditional
expectations, the limit is necessarily $f$.

We now turn to the necessity.
Suppose that for some $1<p<\infty$ and every
$f\in L^p([0,1);X)$, the sequence $(S_Nf)_{N\ge0}$
converges almost everywhere.
By the uniform boundedness principle in $L^0$,
for every $\eta>0$ there exists $C_\eta>0$ such that
\[
 \big|\{x\in[0,1):
       \norm{S_Nf(x)}_X>C_\eta\norm f_{L^p(X)}\}\big|
 \le\eta
\]
for all $N\ge0$ and $f\in L^p([0,1);X)$. Recall from \cite{Wenzel} that
\[
 w_NS_N(w_Nf)
 =\sum_{k\ge0}N_k(\mathbb E_{k+1}-\mathbb  E_k)f,
 \qquad
 N=\sum_{k\ge0}N_k2^k,\quad N_k\in\{0,1\}.
\]
Since $|w_N|=1$, we deduce that
\begin{equation}\label{eq:umd-necessity-measure}
 \left|\left\{x\in[0,1):
 \left\|\sum_{k=1}^m
 \varepsilon_k(\mathbb E_k-\mathbb E_{k-1})f(x)\right\|_X
 >C_\eta\norm f_{L^p(X)}
 \right\}\right|
 \le\eta
\end{equation}
for every $m\ge1$, every choice of
$(\varepsilon_k)_{k=1}^m\in\{0,1\}^m$
and every $f\in L^p([0,1);X)$.

Fix $m\ge1$, an $X$-valued dyadic step function
$f=\mathbb  E_mf$ and $(\varepsilon_k)_{k=1}^m\in\{0,1\}^m$.
Set
\[
 g_j=\sum_{k=1}^j\varepsilon_k(\mathbb E_k-\mathbb E_{k-1})f,
 \qquad 0\le j\le m,
\]
and
\[
 g^*=\max_{0\le j\le m}\norm{g_j}_X,
 \qquad
 F^*=\max_{0\le j\le m}(\mathbb E_j\norm f_X^p)^{1/p}.
\]
We first establish the maximal version of
\eqref{eq:umd-necessity-measure}, namely,
\begin{equation}\label{eq:umd-necessity-maximal}
 |\{g^*>C_\eta\norm f_{L^p(X)}\}|
 \le\eta.
\end{equation}
To this end, consider the stopping time
\[
 \tau=\inf\{0\le j\le m:
       \norm{g_j}_X>C_\eta\norm f_{L^p(X)}\},
\]
with the convention that $\tau=m$ if the set is empty.
The stopped martingale satisfies
\[
 (\mathbb E_k-\mathbb E_{k-1})(\mathbb E_\tau f)
 =\one_{\{\tau\ge k\}}(\mathbb E_k-\mathbb E_{k-1})f,
 \qquad 1\le k\le m.
\]
Consequently,
\[
 g_\tau=\sum_{k=1}^m
 \varepsilon_k(\mathbb E_k-\mathbb E_{k-1})(\mathbb E_\tau f).
\]
Since $\norm{\mathbb E_\tau f}_{L^p(X)}
\le\norm f_{L^p(X)}$, applying
\eqref{eq:umd-necessity-measure} to $\mathbb E_\tau f$ yields
\[
\begin{aligned}
 |\{g^*>C_\eta\norm f_{L^p(X)}\}|
 &=|\{\norm{g_\tau}_X>C_\eta\norm f_{L^p(X)}\}|\\
 &\le
 |\{\norm{g_\tau}_X>C_\eta\norm{\mathbb E_\tau f}_{L^p(X)}\}|\le\eta,
\end{aligned}
\]
which proves \eqref{eq:umd-necessity-maximal}.

We next show that for every $\lambda>0$ and
$0<\delta\le(C_\eta+2)^{-1}$ with $\delta<1/2$,
\begin{equation}\label{eq:umd-necessity-good-lambda}
 |\{g^*>2\lambda,\ F^*\le\delta\lambda\}|
 \le\eta|\{g^*>\lambda\}|.
\end{equation}
Let $\mathcal J$ be the family of pairwise disjoint
dyadic intervals associated with the stopping time
$\inf\{0\le k\le m:\norm{g_k}_X>\lambda\}$.
Then
\[
 \{g^*>\lambda\}
 =\mathop{\dot\bigcup}_{J\in\mathcal J}J.
\]
Fix $J\in\mathcal J$ with $|J|=2^{-k}$ and suppose
that $J\cap\{F^*\le\delta\lambda\}\ne\varnothing$.
By the definition of $F^*$ and the fact that
$\mathbb E_k\norm f_X^p$ and $\mathbb E_{k-1}\norm f_X^p$
are constant on $J$, we have
\[
 \frac1{|J|}\int_J\norm{f(x)}_X^p\,\dd x
 \le(\delta\lambda)^p
\]
and
\[
 \norm{\mathbb E_kf}_X,\ \norm{\mathbb E_{k-1}f}_X
 \le\delta\lambda
 \qquad\text{on}\quad J.
\]
It follows that
\[
\begin{aligned}
 \norm{g_k}_X
 &\le\norm{g_{k-1}}_X
       +\norm{(\mathbb E_k-\mathbb E_{k-1})f}_X \le(1+2\delta)\lambda
 \qquad\text{on}\quad J.
\end{aligned}
\]
Since $\norm{g_j}_X\le\lambda$ on $J$ for $j<k$,
we obtain
\[
 J\cap\{g^*>2\lambda\}
 \subseteq
 \left\{x\in J:
 \max_{k\le j\le m}\norm{g_j(x)-g_k(x)}_X
 >(1-2\delta)\lambda\right\}.
\]
On the other hand, for $k\le j\le m$,
\[
 g_j-g_k
 =\sum_{\ell=k+1}^j
 \varepsilon_\ell(\mathbb E_\ell-\mathbb E_{\ell-1})f.
\]
Rescaling $J$ to $[0,1)$ and applying
\eqref{eq:umd-necessity-maximal}, we obtain
\[
 \left|\left\{x\in J:
 \max_{k\le j\le m}\norm{g_j(x)-g_k(x)}_X
 >C_\eta
 \left(\frac1{|J|}\int_J\norm{f}_X^p\,\dd x
 \right)^{1/p}
 \right\}\right|
 \le\eta|J|.
\]
Since $C_\eta\delta\le1-2\delta$, the preceding
estimates imply
\[
 |J\cap\{g^*>2\lambda\}|\le\eta|J|.
\]
Summing over all $J\in\mathcal J$ which meet
$\{F^*\le\delta\lambda\}$, we arrive at
\[
\begin{aligned}
 |\{g^*>2\lambda,\ F^*\le\delta\lambda\}|
 &\le\eta
 \sum_{\substack{J\in\mathcal J\\
 J\cap\{F^*\le\delta\lambda\}\ne\varnothing}}|J|\le\eta|\{g^*>\lambda\}|.
\end{aligned}
\]
This proves \eqref{eq:umd-necessity-good-lambda}.

Now fix $p<s<\infty$ and choose $\eta<2^{-s}$.
By \eqref{eq:umd-necessity-good-lambda},
\[
 |\{g^*>2\lambda\}|
 \le\eta|\{g^*>\lambda\}|
     +|\{F^*>\delta\lambda\}|.
\]
Multiplying both sides by $s\lambda^{s-1}$
and integrating over $(0,\infty)$, we obtain
\[
 (2^{-s}-\eta)\norm{g^*}_{L^s}^s
 \le\delta^{-s}\norm{F^*}_{L^s}^s.
\]
Thus
\[
 \norm{g^*}_{L^s}
 \le C_{p,s,X}\norm{F^*}_{L^s}.
\]
Since $s/p>1$, it follows from Doob's maximal inequality that
\[
 \norm{F^*}_{L^s}
 =
 \left\|\max_{0\le j\le m}\mathbb E_j\norm f_X^p
 \right\|_{L^{s/p}}^{1/p}
 \le C_{p,s}\norm f_{L^s(X)}.
\]
Combining the above two estimates, we conclude that
\[
 \left\|\sum_{k=1}^m
 \varepsilon_k(\mathbb E_k-\mathbb E_{k-1})f
 \right\|_{L^s(X)}
 \le C_{p,s,X}\norm f_{L^s(X)}
\]
uniformly over $m$ and all choices of
$(\varepsilon_k)_{k=1}^m\in\{0,1\}^m$.

For arbitrary signs
$(\epsilon_k)_{k=1}^m\in\{-1,1\}^m$, writing
$\epsilon_k=2\varepsilon_k-1$ with
$\varepsilon_k\in\{0,1\}$, we obtain
\[
\begin{aligned}
 \left\|\sum_{k=1}^m
 \epsilon_k(\mathbb E_k-\mathbb E_{k-1})f\right\|_{L^s(X)}
 &\le 2\left\|\sum_{k=1}^m
 \varepsilon_k(\mathbb E_k-\mathbb E_{k-1})f\right\|_{L^s(X)}
 +\norm{f-\mathbb E_0f}_{L^s(X)}\\
 &\le C_{p,s,X}\norm f_{L^s(X)}.
\end{aligned}
\]
Applying this estimate to $f-\mathbb E_0f$, we further have
\[
 \left\|\sum_{k=1}^m
 \epsilon_k(\mathbb E_k-\mathbb E_{k-1})f\right\|_{L^s(X)}
 \le C_{p,s,X}\norm{f-\mathbb  E_0f}_{L^s(X)}.
\]
By the characterization of the UMD property in terms
of Paley--Walsh martingales
(see \cite[Theorem~4.2.5]{HNVW}), $X$ is UMD.
The proof is complete.
\end{proof}

\medskip

\section{Decomposition of the linearized model}\label{sec:trees}\label{subsec:haar-multipliers}\label{subsec:coefficient-selection}\label{sec:selection}\label{subsec:joint-levels}

By the reductions established in Section~\ref{sec:reduction}, the proof of
Theorem~1.1 is reduced to establishing the  $L^p$
bound for the  linearization operators $L_{\PP}$ in
Theorem~4.4. The purpose of this section is to decompose
each fixed linearization into a finite sum of operators
associated with trees, with suitable estimates on their
size, density and top intervals. This decomposition will
be used in the next section to prove Theorem~4.4.

Fix $1<r<\infty$, a finite bitile family
$\PP\subseteq\mathcal B$  and a linearization $L_{\PP}$ as in
Subsection~\ref{subsec:linearization}. Thus, all linearization parameters (namely, the signs
$(\epsilon_P)_{P\in\PP}$, the integer $K$, the measurable map
$\mathbf N=(N_0,\ldots,N_K)$  and the strongly measurable functions
$a_j^*$, $1\le j\le K$) are fixed and satisfy \eqref{eq:directions}.
For a finite subfamily $\mathbb A\subseteq\PP$ and
$g\in L^{r'}(\R_+)$, define its density by
\begin{equation}\label{eq:density}
\begin{split}
\dens(\mathbb A;g)
:=\sup_{P\in\mathbb A}\ \sup_{\substack{Q\in\mathcal B\\P\le Q}}
\left(\frac1{|I_Q|}\int_{I_Q}|g(x)|^{r'}
\sum_{j=1}^K\norm{a_j^*(x)}_{X^*}^{r'}
\one_{\omega_Q}(N_j(x))\dd x\right)^{1/r'}.
\end{split}
\end{equation}
We set $\dens(\varnothing;g):=0$. Clearly, this definition gives
$\dens(\mathbb A';g)\le\dens(\mathbb A;g)$ whenever
$\mathbb A'\subseteq\mathbb A$. Moreover, if $|g|\le\one_E$ for a
measurable set $E\subseteq\R_+$, then by \eqref{eq:directions}, we have
\[
|g(x)|^{r'}\sum_{j=1}^K\norm{a_j^*(x)}_{X^*}^{r'}
\one_{\omega_Q}(N_j(x))\le\one_E(x),
\]
and hence $\dens(\mathbb A;g)\le1$. For $2\le q<\infty$ and $f\in L^1_{\mathrm{loc}}(\R_+;X)$, define
\begin{equation}\label{eq:size}
\size_q(\mathbb A;f)
:=\sup_{\substack{\varnothing\ne\VV\subseteq\mathbb A,\ T\in\mathcal B\\
P\le_u T\ \text{for every }P\in\VV}}
|I_T|^{-1/q}
\left\|\sum_{P\in\VV}\ip f{w_{P_d}}w_{P_d}\right\|_{L^q(\R_+;X)}.
\end{equation}
Here the supremum is taken over all nonempty subfamilies contained in an
up-tree with top $T$. We also put $\size_q(\varnothing;f):=0$.
For $\mathbb A'\subseteq\mathbb A$, we have
$\size_q(\mathbb A';f)\le\size_q(\mathbb A;f)$. Taking $\VV=\{P\}$ and $T=P$
in \eqref{eq:size}, and using
$|w_{P_d}|=|I_P|^{-1/2}\one_{I_P}$, we obtain
\begin{equation}\label{eq:singleton}
|I_P|^{-1/2}\norm{\ip f{w_{P_d}}}_X
=|I_P|^{-1/q}
\norm{\ip f{w_{P_d}}w_{P_d}}_{L^q(X)}
\le\size_q(\mathbb A;f),\quad   P\in\mathbb A.
\end{equation}
For every subfamily $\mathbb A\subseteq\PP$, we denote by
$L_{\mathbb A}$ the operator obtained from $L_{\PP}$ by
restricting the sum in \eqref{eq:linearization} to
$P\in\mathbb A$, with the same signs and linearization
parameters.  Our  main result of this section reads as follows.

\begin{theorem}\label{prop:joint-levels}
	Let $2\le q<r<\infty$, and suppose that $X$ is a UMD
	space satisfying $(\mathrm W_q)$ with constant $A_q$.
	Let $\PP\subseteq\mathcal B$ be a finite bitile family
	and fix a linearization $L_{\PP}$ as in
	Subsection~\ref{subsec:linearization}. Let
	$E,F\subseteq\R_+$ be measurable sets with
	$0<|E|,|F|<\infty$. Suppose that
	$f\in L^q(\R_+;X)$ and $g\in L^{r'}(\R_+)$ satisfy
	\[
	\|f(x)\|_X\le \one_F(x),
	\qquad
	|g(x)|\le \one_E(x)
	\]
	for almost every $x\in\R_+$. Then there exist a
	(possibly empty) subfamily $\PP_{\mathrm{rem}}\subseteq\PP$
	and finitely many nonempty trees $\mathbb T_{n,j}$ such
	that the following  hold.
	
	\begin{enumerate}[\rm (i)]
		\item
		The family $\PP$ admits the disjoint decomposition
		\[
		\PP
		=
		\PP_{\mathrm{rem}}
		\mathbin{\dot\cup}
		\mathop{\dot\bigcup}_{n\in\Z}
		\mathop{\dot\bigcup}_{j}\mathbb T_{n,j}.
		\]
		Each $\mathbb T_{n,j}$ is a tree with top
		$T_{n,j}\in\mathcal B$ satisfying
		\[
		\mathbb T_{n,j}\subseteq\PP,
		\qquad
		P\le T_{n,j}
		\quad(P\in\mathbb T_{n,j}).
		\]
		
		\item
		The remainder satisfies
		\[
		\langle f,w_{P_d}\rangle=0
		\qquad(P\in\PP_{\mathrm{rem}}).
		\]
		
		\item
		For every $n\in\Z$ and every tree $\mathbb T_{n,j}$,
		\begin{align}
			\operatorname{size}_q(\mathbb T_{n,j};f)
			&\le C\min\left\{
			\operatorname{size}_q(\PP;f),
			2^{n/q}|F|^{1/q}
			\right\},
			\label{eq:level-size}\\
			\operatorname{dens}(\mathbb T_{n,j};g)
			&\le C\min\left\{
			\operatorname{dens}(\PP;g),
			2^{n/r'}|E|^{1/r'}
			\right\}.
			\label{eq:level-density}
		\end{align}
		
		\item
		The top intervals satisfy
		\begin{align}
			\sum_j |I_{T_{n,j}}|
			&\le C2^{-n},
			\qquad n\in\Z,
			\label{eq:level-total}\\
			\sum_{j:\,I_{T_{n,j}}\subseteq J}
			|I_{T_{n,j}}|
			&\le C2^{-n}
			\left(
			\frac{|F\cap J|}{|F|}
			+
			\frac{|E\cap J|}{|E|}
			\right),
			\qquad J\in\mathcal D,\ n\in\Z.
			\label{eq:level-local}
		\end{align}
		The tree decomposition of $\PP$ in {\rm (i)} is chosen
		independently of $J$, so that
		\eqref{eq:level-local} holds simultaneously for all
		$J\in\mathcal D$. The top intervals need not be
		pairwise disjoint; their lengths are counted with
		multiplicity in both estimates.
		
		\item
		With the same signs and linearization parameters
		on every subfamily, we have
		\[
		L_{\PP_{\mathrm{rem}}}f=0
		\qquad\text{and}\qquad
		L_{\PP}f
		=
		\sum_{n,j}L_{\mathbb T_{n,j}}f.
		\]
	\end{enumerate}
	
	The constant $C$ depends only on
	$q,r,X$ and $A_q$.
\end{theorem}

We prove Theorem~\ref{prop:joint-levels} by first selecting trees with
large density and then selecting trees with large size. These two selection
results are stated separately below.

\begin{proposition}[Density selection]\label{lem:density-selection}
	Fix $r$ and the linearization parameters. Let
	$\PP\subseteq\mathcal B$ be a finite bitile family,
	$g\in L^{r'}(\R_+)$  and $\Delta>0$.
	Then $\PP$ admits a disjoint decomposition
	\[
	\PP
	=
	\PP_{\mathrm{sparse}}
	\mathbin{\dot\cup}
	\mathop{\dot\bigcup}_j\TT_j,
	\]
	where $\PP_{\mathrm{sparse}}\subseteq\PP$ is possibly empty
	and each $\TT_j$ is a nonempty tree with top
	$T_j\in\mathcal B$. In particular,
$P\le T_j$ for all  $P\in\TT_j$.
	Moreover,
	\[
	\dens(\PP_{\mathrm{sparse}};g)\le\Delta,
	\]
	and the tops are pairwise disjoint as bitiles, with
	\begin{equation}\label{eq:local-density-pack}
		\sum_{j:I_{T_j}\subseteq J}|I_{T_j}|
		\le\Delta^{-r'}\int_J|g|^{r'}\dd x,
		\qquad J\in\mathcal D.
	\end{equation}
\end{proposition}

%
%

\begin{proof}
For each $P\in\PP$ with $\dens(\{P\};g)>\Delta$, choose a bitile
$P'\ge P$ such that
\[
\Delta^{r'}|I_{P'}|
<\int_{I_{P'}}|g(x)|^{r'}
\sum_{j=1}^K\norm{a_j^*(x)}_{X^*}^{r'}
\one_{\omega_{P'}}(N_j(x))\dd x.
\]
The chosen bitiles form a finite family. After removing repetitions,
let $T_1,\ldots,T_m$ be its maximal members with respect to $\le$.
Assign the members of $\PP$ to these tops in order by setting
\[
\begin{aligned}
&\TT_i=\left\{P\in\PP\setminus\bigcup_{k<i}\TT_k:P\le T_i\right\},
\qquad 1\le i\le m,\\
&\PP_{\mathrm{sparse}}=\PP\setminus\bigcup_{i=1}^m\TT_i.
\end{aligned}
\]
If $\dens(\{P\};g)>\Delta$, its chosen bitile $P'$ satisfies
$P\le P'\le T_i$ for some maximal member $T_i$. The assignment
therefore places $P$ in one of the trees, and
\[
\{P\in\PP:\dens(\{P\};g)>\Delta\}
\subseteq\bigcup_{i=1}^m\TT_i,
\]
which, by the definition of density, gives
\[
\dens(\PP_{\mathrm{sparse}};g)
=\sup_{P\in\PP_{\mathrm{sparse}}}\dens(\{P\};g)
\le\Delta.
\]
The same conclusion is immediate if $\PP_{\mathrm{sparse}}=\varnothing$.

By maximality, we have $T_i\not\le T_k$ and $T_k\not\le T_i$ for $i\ne k$.
Since distinct bitiles that intersect are comparable, it follows that
\[
T_i\cap T_k=\varnothing\qquad(i\ne k).
\]
In particular, for every $J\in\mathcal D$, every $x\in\R_+$, and
$1\le j\le K$,
\[
\sum_{i:I_{T_i}\subseteq J}
\one_{I_{T_i}}(x)\one_{\omega_{T_i}}(N_j(x))
\le\one_J(x).
\]
Using the inequalities defining the selected tops, we obtain
\begin{align*}
\Delta^{r'}\sum_{i:I_{T_i}\subseteq J}|I_{T_i}|
&\le\sum_{i:I_{T_i}\subseteq J}
\int_{I_{T_i}}|g(x)|^{r'}
\sum_{j=1}^K\norm{a_j^*(x)}_{X^*}^{r'}
\one_{\omega_{T_i}}(N_j(x))\dd x\\
&\le\int_J|g(x)|^{r'}
\sum_{j=1}^K\norm{a_j^*(x)}_{X^*}^{r'}\dd x\le\int_J|g(x)|^{r'}\dd x,
\end{align*}
where the last inequality follows from \eqref{eq:directions}. Discarding
empty trees does not increase the sum of the top lengths. The remaining
trees and tops therefore satisfy the asserted estimates. The proof is complete.
\end{proof}

\begin{proposition}[Size selection]\label{lem:size-selection}
	Let $2\le q<\infty$. Suppose that $X$ is UMD
	satisfying $(\mathrm W_q)$ with constant $A_q$.
	Let $\PP\subseteq\mathcal B$ be a finite bitile family,
	$f\in L^q(\R_+;X)$ and $\sigma>0$.
	Then there exist a (possibly empty) subfamily
	$\PP_{\mathrm{small}}\subseteq\PP$ and finitely many
	nonempty trees $\TT_j$ with tops
	$T_j\in\PP^{\mathrm{cvx}}$ such that
	\[
	\PP
	=
	\PP_{\mathrm{small}}
	\mathbin{\dot\cup}
	\mathop{\dot\bigcup}_j\TT_j,
	\]
	where each $\TT_j\subseteq\PP$ satisfies $
	P\le T_j$ for all $P\in\TT_j.$
	Moreover,
\begin{equation}\label{eq:size-estimate}	\size_q(\PP_{\mathrm{small}};f)\le\sigma,
\end{equation}
	and the top intervals satisfy
	\begin{equation}\label{eq:local-size-pack}.
		\sum_{j:I_{T_j}\subseteq J}|I_{T_j}|
		\le (2\beta_{q,X})^q A_q^q\sigma^{-q}
		\int_J\norm{f(x)}_X^q\dd x,\qquad J\in\mathcal D,
	\end{equation}
where $\beta_{q,X}=4\max\{1,C_{q,X}\}$ and
$C_{q,X}$ is the UMD constant.
\end{proposition}

The proof of Proposition~\ref{lem:size-selection} relies on two preparatory lemmas. The first provides a bounded function with the same up- and down-packet coefficients on a finite convex bitile family.

\begin{lemma}\label{lem:reconstruction}
Let $\Tiles\subseteq\mathcal B$ be a finite convex
bitile family and let $T\in\mathcal B$.
Suppose that $f\in L^1(I_T;X)$ and $M\ge0$ satisfy
\[
    |I_P|^{-1/2}\norm{\ip f{w_{P_\nu}}}_X\le M,
    \qquad P\in\Tiles,\ P\le T,\ \nu\in\{d,u\}.
\]
Then there exists an $X$-valued step function $g_T$
supported on $I_T$ such that
\begin{equation}\label{eq:reconstruction-0}
\norm{g_T}_{L^\infty(X)}\le2M
\end{equation}
and
\begin{equation}\label{eq:reconstruction}
\begin{aligned}
    \ip{g_T}{w_{P_\nu}}
    &=\ip f{w_{P_\nu}},
    \qquad P\in\Tiles,\ P\le T,\ \nu\in\{d,u\}.
\end{aligned}
\end{equation}
\end{lemma}

\begin{proof}
Put
\[
\mathcal I_T=\{I_P:P\in\Tiles,\ P\le T\},
\qquad
\mathcal K_T
=\{K\in\mathcal D:K^{(1)}\in\mathcal I_T,\
K\notin\mathcal I_T\}.
\]
If $\mathcal I_T=\varnothing$, then one may take $g_T=0$. Hence, we assume from now on $\mathcal I_T\ne\varnothing$. We first show that $\mathcal I_T$ is convex.
For $I,J\in\mathcal I_T$ and a dyadic interval
$L$ with $I\subseteq L\subseteq J$, let
$P_I,P_L,P_J$ be the unique bitiles below $T$
with time intervals $I,L,J$, respectively.
Then
\[
P_I\le P_L\le P_J.
\]
Since $P_I,P_J\in\Tiles$ and $\Tiles$ is convex,
we have $P_L\in\Tiles$, and hence $L\in\mathcal I_T$. We now claim that the intervals in $\mathcal K_T$ are pairwise disjoint. Indeed, if two distinct intervals $K,L\in\mathcal K_T$
intersect, we may assume $K\subsetneq L$. Then
$
K^{(1)}\subseteq L\subsetneq L^{(1)}$ and $
K^{(1)},L^{(1)}\in\mathcal I_T.
$
By convexity, $L\in\mathcal I_T$, a contradiction. Moreover, for every $I\in\mathcal I_T$,
\begin{equation}\label{eq:terminal-partition}
I_\pm
=\mathop{\dot\bigcup}_{\substack{K\in\mathcal K_T\\
K\subseteq I_\pm}}K.
\end{equation}
To verify this, fix $x\in I_\pm$ and choose the
smallest $J\in\mathcal I_T$ such that
$x\in J\subseteq I$. Let $K$ be the child of $J$
containing $x$. By minimality,
$
K^{(1)}=J\in\mathcal I_T$ and $K\notin\mathcal I_T.$
Since $x\in K\cap I_\pm$ and
$|K|\le |I_\pm|$,   dyadic nesting implies
$K\subseteq I_\pm$. Thus
\eqref{eq:terminal-partition} follows from the
pairwise disjointness of $\mathcal K_T$.

Now take $m_T=|I_T|^{1/2}w_{T_u}$ and set
$f_T=m_Tf$.
Define
\[
f_T^0=\sum_{K\in\mathcal K_T}
\langle f_T\rangle_K\one_K,
\qquad
g_T=m_Tf_T^0.
\]
By \eqref{eq:two-components} and the hypothesis,
\[
\norm{\langle f_T\rangle_I}_X\le M,
\qquad
|I|^{-1/2}\norm{\ip{f_T}{h_I}}_X\le M,\qquad I\in\mathcal I_T.
\]
Consequently,
\[
\langle f_T\rangle_{I_\pm}
=\langle f_T\rangle_I
\mp |I|^{-1/2}\ip{f_T}{h_I},
\qquad
\norm{\langle f_T\rangle_{I_\pm}}_X\le2M.
\]
Since every $K\in\mathcal K_T$ is a child of
some interval in $\mathcal I_T$, we have
\[
\norm{g_T}_{L^\infty(X)}
=\norm{f_T^0}_{L^\infty(X)}
\le2M,
\]
which gives the estimate \eqref{eq:reconstruction-0}.
Finally, by \eqref{eq:terminal-partition},
for every $I\in\mathcal I_T$,
\[
\int_{I_\pm}f_T^0
=\sum_{\substack{K\in\mathcal K_T\\
K\subseteq I_\pm}}
|K|\langle f_T\rangle_K
=\int_{I_\pm}f_T.
\]
Consequently,
\[
\langle f_T^0\rangle_I=\langle f_T\rangle_I,
\qquad
\ip{f_T^0}{h_I}=\ip{f_T}{h_I}.
\]
By \eqref{eq:two-components}, these identities imply
\[
\ip{f_T^0}{m_Tw_{P_\nu}}
=\ip{f_T}{m_Tw_{P_\nu}},
\qquad P\in\Tiles,\ P\le T,\ \nu\in\{d,u\}.
\]
Since $m_Tg_T=f_T^0$, $f_T=m_Tf$  and
$m_T^2=\one_{I_T}$, we further get
\[
\begin{aligned}
\ip{g_T}{w_{P_\nu}}
&=\ip{m_Tg_T}{m_Tw_{P_\nu}}\\
&=\ip{f_T^0}{m_Tw_{P_\nu}}\\
&=\ip{f_T}{m_Tw_{P_\nu}}
=\ip f{w_{P_\nu}},
\end{aligned}
\]
which gives the desired estimate \eqref{eq:reconstruction}.
\end{proof}

The following lemma combines the bounded reconstruction
of Lemma~\ref{lem:reconstruction} with the UMD
property to control the size
by the packet coefficients on the convex hull.

\begin{lemma}\label{prop:size-bound}
Let $X$ be a complex UMD space, $2\le q<\infty$,
and $\PP\subseteq\mathcal B$ a finite bitile family.
Suppose that $f\in L^1_{\mathrm{loc}}(\R_+;X)$
and $M\ge0$ satisfy
\[
|I_R|^{-1/2}\norm{\ip f{w_{R_\nu}}}_X\le M,
\qquad R\in\PP^{\mathrm{cvx}},\ \nu\in\{d,u\}.
\]
Then
\[
\size_q(\PP;f)\le2\beta_{q,X}M,
\]
where $\beta_{q,X}$ is the same as in Proposition \ref{lem:size-selection}.
In particular, if
$\norm{f(\cdot)}_X\le\one_F$ almost everywhere
for some measurable $F\subseteq\R_+$, then
\[
\size_q(\PP;f)\le2\beta_{q,X}.
\]
\end{lemma}

\begin{proof}
Note that $\PP^{\mathrm{cvx}}$ is finite and convex,
since between two comparable bitiles there is at most one intermediate bitile at each dyadic time scale. Let $\VV\subseteq\PP$ be a nonempty up-tree with top $T$.
By hypothesis,
\[
|I_P|^{-1/2}\norm{\ip f{w_{P_\nu}}}_X\le M,
\qquad
P\in\PP^{\mathrm{cvx}},\ P\le T,\ \nu\in\{d,u\}.
\]
Applying Lemma~\ref{lem:reconstruction}, there exists an
$X$-valued step function $g_T$, supported on $I_T$,
such that
\[
\norm{g_T}_{L^\infty(X)}\le2M,
\qquad
\ip{g_T}{w_{P_\nu}}=\ip f{w_{P_\nu}},
\quad
P\in\PP^{\mathrm{cvx}},\ P\le T,\ \nu\in\{d,u\}.
\]
Take $m_T=|I_T|^{1/2}w_{T_u}$.
By \eqref{eq:two-components},
\[
m_Tw_{P_d}=\eta_Ph_{I_P},
\qquad
\eta_P\in\{-1,1\},
\quad P\in\VV.
\]
The coefficient identities above give
\[
\begin{aligned}
m_T\sum_{P\in\VV}\ip f{w_{P_d}}w_{P_d}
&=\sum_{P\in\VV}
\ip{g_T}{w_{P_d}}m_Tw_{P_d}\\
&=\sum_{P\in\VV}
\ip{m_Tg_T}{h_{I_P}}h_{I_P}.
\end{aligned}
\]
Since the Walsh packet sum is supported on $I_T$
and $|m_T|=1$ on $I_T$, we obtain
\[
\begin{aligned}
\left\|
\sum_{P\in\VV}\ip f{w_{P_d}}w_{P_d}
\right\|_{L^q(X)}
&=
\left\|
m_T\sum_{P\in\VV}\ip f{w_{P_d}}w_{P_d}
\right\|_{L^q(X)}\\
&=
\left\|
\sum_{P\in\VV}
\ip{m_Tg_T}{h_{I_P}}h_{I_P}
\right\|_{L^q(X)}.
\end{aligned}
\]
Write $|I_T|=2^{-k_0}$ and choose $k_1\ge k_0$
such that
\[
2^{-k_1}\le |I_P|\le2^{-k_0},
\qquad P\in\VV.
\]
Distinct bitiles below $T$ have distinct time
intervals, which are pairwise disjoint at each
fixed scale. Consequently,
\[
\begin{aligned}
&\sum_{P\in\VV}
\ip{m_Tg_T}{h_{I_P}}h_{I_P}=
\sum_{k=k_0}^{k_1}
\left(
\sum_{\substack{P\in\VV\\|I_P|=2^{-k}}}
\one_{I_P}
\right)
(\E_{k+1}-\E_k)(m_Tg_T).
\end{aligned}
\]
For each $k$, the multiplier in parentheses is
$\mathscr D_k$-measurable and takes values in
$\{0,1\}$. Thus the right-hand side is a
predictable martingale transform. Applying the UMD martingale-transform inequality
 and the contractivity of
conditional expectations, we obtain
\[
\begin{aligned}
\left\|
\sum_{P\in\VV}
\ip f{w_{P_d}}w_{P_d}
\right\|_{L^q(X)}&=
\left\|
\sum_{P\in\VV}
\ip{m_Tg_T}{h_{I_P}}h_{I_P}
\right\|_{L^q(X)}\\
&\quad\le
C_{q,X}
\norm{(\E_{k_1+1}-\E_{k_0})(m_Tg_T)}_{L^q(X)}\\
&\quad\le
2C_{q,X}\norm{m_Tg_T}_{L^q(X)}.
\end{aligned}
\]
Note that
$g_T$ is supported on $I_T$ and
$\norm{g_T}_{L^\infty(X)}\le2M$.
Since $|m_T|=1$ on $I_T$, we conclude
\[
\norm{m_Tg_T}_{L^q(X)}
=\norm{g_T}_{L^q(X)}
\le2M|I_T|^{1/q}.
\]
Combining the above two estimates, we get
\begin{equation}\label{eq:masked-tree}
\left\|
\sum_{P\in\VV}
\ip f{w_{P_d}}w_{P_d}
\right\|_{L^q(X)}
\le2\beta_{q,X}M|I_T|^{1/q}.
\end{equation}
Taking the supremum over all up-trees
$\VV\subseteq\PP$ and their tops, we arrive at
\[
\size_q(\PP;f)\le2\beta_{q,X}M.
\]

Finally, suppose that
$\norm{f(\cdot)}_X\le\one_F$ almost everywhere.
For every $R\in\PP^{\mathrm{cvx}}$ and
$\nu\in\{d,u\}$,
\[
\begin{aligned}
|I_R|^{-1/2}\norm{\ip f{w_{R_\nu}}}_X
&\le\frac1{|I_R|}
\int_{I_R}\norm{f(x)}_X\dd x\le\frac{|F\cap I_R|}{|I_R|}
\le1.
\end{aligned}
\]
Applying the preceding estimate with $M=1$
completes the proof.
\end{proof}

With the reconstruction lemma and the size estimate
established above, we are now ready to prove
Proposition~\ref{lem:size-selection}.

\begin{proof}[Proof of Proposition~\ref{lem:size-selection}]
Let
\[
\PP_{\max}=\left\{T\in\PP^{\mathrm{cvx}}:
\max_{\nu\in\{d,u\}}|I_T|^{-1/2}
\norm{\ip f{w_{T_\nu}}}_X
>\frac{\sigma}{2\beta_{q,X}}\right\}.
\]
Since $\PP^{\mathrm{cvx}}$ is finite, we may enumerate
the maximal elements of $\PP_{\max}$ with respect to
$\le$ as $T_1,\ldots,T_m$. If $\PP_{\max}=\varnothing$,
we take $m=0$, with the convention that unions over
$1\le i\le m$ are empty. For each $1\le i\le m$, choose
$Q_i\in\{(T_i)_d,(T_i)_u\}$ such that
\[
|I_{T_i}|^{-1/2}\norm{\ip f{w_{Q_i}}}_X
>\frac{\sigma}{2\beta_{q,X}}.
\]
The maximal bitiles $T_i$ are pairwise disjoint,
and hence so are their component tiles $Q_i$. We now construct the required decomposition.
For $1\le i\le m$, define
\[
\TT_i=
\left\{P\in\PP\setminus\bigcup_{k<i}\TT_k:
P\le T_i\right\},
\]
and set
\[
\PP_{\mathrm{small}}
=\PP\setminus\bigcup_{i=1}^m\TT_i.
\]
Then
\[
\PP=\PP_{\mathrm{small}}
\mathbin{\dot\cup}\mathop{\dot\bigcup}_{i=1}^m\TT_i,
\]
where each nonempty $\TT_i$ is a tree with top $T_i$.

We first verify the size estimate \eqref{eq:size-estimate}. Let
$R\in\PP_{\mathrm{small}}^{\mathrm{cvx}}$.
Choose $P,P'\in\PP_{\mathrm{small}}$ with
$P\le R\le P'$. Then $R\in\PP^{\mathrm{cvx}}$.
We claim that $R\notin\PP_{\max}$. Otherwise, by
the maximality of $T_1,\ldots,T_m$, there is some
$i$ such that $R\le T_i$. Consequently,
$P\le T_i$, and since $P$ has not been assigned
to any preceding tree, the definition of $\TT_i$
gives $P\in\TT_i$, contrary to
$P\in\PP_{\mathrm{small}}$. This verifies the claim, and so for any $f\in L^q(\R_+;X)$ we have
\[
|I_R|^{-1/2}\norm{\ip f{w_{R_\nu}}}_X
\le\frac{\sigma}{2\beta_{q,X}},
\qquad \nu\in\{d,u\}.
\]
Applying Lemma~\ref{prop:size-bound} to
$\PP_{\mathrm{small}}$, we obtain
\[
\size_q(\PP_{\mathrm{small}};f)
\le2\beta_{q,X}\frac{\sigma}{2\beta_{q,X}}
=\sigma,
\]
which gives the desired estimate \eqref{eq:size-estimate}.

It remains to prove the local packing estimate \eqref{eq:local-size-pack}.
Fix $J\in\mathcal D$. For every $i$ with
$I_{T_i}\subseteq J$, we have
$
\ip f{w_{Q_i}}
=\ip{f\one_J}{w_{Q_i}},
$
since $w_{Q_i}$ is supported on $I_{T_i}$.
Moreover, using
$|w_{Q_i}|=|I_{T_i}|^{-1/2}\one_{I_{T_i}}$
and the choice of $Q_i$, we obtain
\[
\begin{aligned}
\norm{\ip{f\one_J}{w_{Q_i}}w_{Q_i}}_{L^q(X)}^q
&=\norm{\ip f{w_{Q_i}}}_X^q
|I_{T_i}|^{1-q/2}>\left(\frac{\sigma}{2\beta_{q,X}}\right)^q
|I_{T_i}|.
\end{aligned}
\]
Applying $(\mathrm W_q)$ to $f\one_J$ and the
pairwise disjoint tiles $Q_i$ with $I_{T_i}\subseteq J$,
we get
\[
\begin{aligned}
\left(\frac{\sigma}{2\beta_{q,X}}\right)^q
\sum_{i:I_{T_i}\subseteq J}|I_{T_i}|
&\le
\sum_{i:I_{T_i}\subseteq J}
\norm{\ip{f\one_J}{w_{Q_i}}w_{Q_i}}_{L^q(X)}^q\\
&\le A_q^q\norm{f\one_J}_{L^q(X)}^q=A_q^q\int_J\norm{f(x)}_X^q\dd x.
\end{aligned}
\]
Thus the required local packing estimate \eqref{eq:local-size-pack} holds
for all the selected tops $T_i$.

Finally, discarding the empty trees does not affect the
decomposition or the preceding estimates.
This completes the proof.
\end{proof}

We now prove  Theorem~\ref{prop:joint-levels} by combining the density and size
selection propositions to obtain the required decomposition.

\begin{proof}[Proof of Theorem~\ref{prop:joint-levels}]
	If $\ip f{w_{P_d}}=0$ for every $P\in\PP$, take
	$\PP_{\mathrm{rem}}=\PP$ and no trees. Then
	$L_{\PP}f=0$ by \eqref{eq:linearization} and all the
	assertions follow. We may therefore assume that at least
	one down-packet coefficient is nonzero. 	By the support assumptions on $f$ and $g$, Lemma~\ref{prop:size-bound}, \eqref{eq:density} and \eqref{eq:directions}, we have
	\[
	\size_q(\PP;f)\le2\beta_{q,X},
	\qquad
	\dens(\PP;g)\le1.
	\]
	Set
	\[
	\mu=
	\min_{\substack{P\in\PP\\
			\ip f{w_{P_d}}\ne0}}
	|I_P|^{-1/2}\norm{\ip f{w_{P_d}}}_X.
	\]
	Since $\PP$ is finite and at least one down-packet
	coefficient is nonzero, it follows that $\mu>0$. Now take $
	n=\left\lceil\log_2\frac{\mu^q}{|F|}\right\rceil-1.$
	Then $
	2^{n/q}|F|^{1/q}<\mu.$
	Choose an integer $n_0>n$ sufficiently large so that
	\[
	\size_q(\PP;f)\le2^{n_0/q}|F|^{1/q},
	\qquad
	\dens(\PP;g)\le2^{n_0/r'}|E|^{1/r'}.
	\]
	
We now construct the required decomposition by successive
applications of Propositions~\ref{lem:density-selection}
and~\ref{lem:size-selection}. Set
$\mathcal R_{n_0}=\PP$. Starting from
$\mathcal R_{n_0}$, we proceed recursively in decreasing
order of the levels
$
k=n_0-1,n_0-2,\ldots,n.
$
At each level $k$, we first apply
Proposition~\ref{lem:density-selection} to
$\mathcal R_{k+1}$ with
$
\Delta=2^{k/r'}|E|^{1/r'}.
$
Denote by $m_k$ the number of nonempty trees obtained. Then
\[
\mathcal R_{k+1}
=
\mathcal R_{k+1,\mathrm{sparse}}
\mathbin{\dot\cup}
\mathop{\dot\bigcup}_{j=1}^{m_k}\TT_{k,j}.
\]
Next, apply Proposition~\ref{lem:size-selection} to
$\mathcal R_{k+1,\mathrm{sparse}}$ with
$
\sigma=2^{k/q}|F|^{1/q}.
$
	Let $M_k$ denote the total number of trees obtained
at level $k$ after these two selections. Thus the
size selection produces $M_k-m_k$ additional
nonempty trees. Writing the resulting remainder as
$\mathcal R_k$ and indexing these trees after those
obtained in the first step, we have
\[
\mathcal R_{k+1,\mathrm{sparse}}
=
\mathcal R_k
\mathbin{\dot\cup}
\mathop{\dot\bigcup}_{j=m_k+1}^{M_k}\TT_{k,j}.
\]
Consequently,
\[
\mathcal R_{k+1}
=
\mathcal R_k
\mathbin{\dot\cup}
\mathop{\dot\bigcup}_{j=1}^{M_k}\TT_{k,j},
\qquad
\mathcal R_k
\subseteq\mathcal R_{k+1,\mathrm{sparse}}
\subseteq\mathcal R_{k+1}.
\]
After completing the $n_0-n$ steps of the recursion,
we set $\PP_{\mathrm{rem}}=\mathcal R_n$. Combining the
successive disjoint decompositions, we obtain
\[
\PP
=
\PP_{\mathrm{rem}}
\mathbin{\dot\cup}
\mathop{\dot\bigcup}_{k=n}^{n_0-1}
\mathop{\dot\bigcup}_{j=1}^{M_k}\TT_{k,j}.
\]

We now verify the required properties. The decomposition and tree properties in (i) follow
directly from the construction. To prove (ii), suppose, contrary to the assertion, that
there exists $P\in\PP_{\mathrm{rem}}$ such that
$\ip f{w_{P_d}}\ne0$. By the choice of $n$,
\eqref{eq:singleton} and Proposition~\ref{lem:size-selection}
at level $n$, we have
\[
\begin{aligned}
	2^{n/q}|F|^{1/q}
	&<
	|I_P|^{-1/2}\norm{\ip f{w_{P_d}}}_X\le\size_q(\mathcal R_n;f)
	\le2^{n/q}|F|^{1/q},
\end{aligned}
\]
which is a contradiction. Hence
\[
\ip f{w_{P_d}}=0
\qquad(P\in\PP_{\mathrm{rem}}).
\]
We next establish the size and density estimates in
(iii). By the two selection estimates and the inclusions
$
\mathcal R_k
\subseteq\mathcal R_{k+1,\mathrm{sparse}}
\subseteq\mathcal R_{k+1},
$
we have, for $n\le k<n_0$,
\[
\size_q(\mathcal R_k;f)
\le2^{k/q}|F|^{1/q},
\qquad
\dens(\mathcal R_k;g)
\le2^{k/r'}|E|^{1/r'}.
\]
The corresponding bounds for $\mathcal R_{n_0}=\PP$
follow from the choice of $n_0$.
Since $\TT_{k,j}\subseteq\mathcal R_{k+1}$, it follows from the
monotonicity of size and density that
\[
\begin{aligned}
	\size_q(\TT_{k,j};f)
	&\le
	\min\left\{
	\size_q(\PP;f),
	2^{(k+1)/q}|F|^{1/q}
	\right\}\\
	&\le
	2^{1/q}\min\left\{
	\size_q(\PP;f),
	2^{k/q}|F|^{1/q}
	\right\},\\
	\dens(\TT_{k,j};g)
	&\le
	\min\left\{
	\dens(\PP;g),
	2^{(k+1)/r'}|E|^{1/r'}
	\right\}\\
	&\le
	2^{1/r'}\min\left\{
	\dens(\PP;g),
	2^{k/r'}|E|^{1/r'}
	\right\},
\end{aligned}
\]
which gives the  size and density estimates in (iii). We turn to establish the packing estimates in (iv).
Fix $J\in\mathcal D$ and $n\le k<n_0$.
The first $m_k$ trees were selected by the density
selection proposition, while the remaining trees
were selected by the size selection proposition.
Consequently, by \eqref{eq:local-density-pack} and
\eqref{eq:local-size-pack},
\[
\begin{aligned}
	\sum_{j:I_{T_{k,j}}\subseteq J}|I_{T_{k,j}}|
	&=
	\sum_{\substack{1\le j\le m_k\\
			I_{T_{k,j}}\subseteq J}}|I_{T_{k,j}}|
	+
	\sum_{\substack{m_k<j\le M_k\\
			I_{T_{k,j}}\subseteq J}}|I_{T_{k,j}}|\\
	&\le
	\frac{2^{-k}}{|E|}
	\int_J|g(x)|^{r'}\dd x\\
	&\quad+
	(2\beta_{q,X})^qA_q^q
	\frac{2^{-k}}{|F|}
	\int_J\norm{f(x)}_X^q\dd x\\
	&\le
	C2^{-k}\left(
	\frac{|E\cap J|}{|E|}
	+
	\frac{|F\cap J|}{|F|}
	\right).
\end{aligned}
\]
The construction is independent of $J$, so this estimate
holds simultaneously for all $J\in\mathcal D$.
This proves the local packing estimate
\eqref{eq:level-local}. To obtain the global packing estimate  \eqref{eq:level-total}, choose
$J=[0,2^m)$ large enough to contain all the top intervals
at level $k$. Then
\[
\sum_{j=1}^{M_k}|I_{T_{k,j}}|
\le2C2^{-k},
\]
which gives \eqref{eq:level-total} after changing
the constant. The estimate is trivial if the level
contains no trees. Thus, (iv) is verified. Finally, by (ii) and \eqref{eq:linearization},
$
L_{\PP_{\mathrm{rem}}}f=0.
$
All the subfamilies retain the same signs and
linearization parameters $K$, $\mathbf N$, and $a_j^*$.
Thus, by the disjoint decomposition of $\PP$ and
\eqref{eq:linearization},
\[
\begin{aligned}
	L_{\PP}f
	&=L_{\PP_{\mathrm{rem}}}f
	+\sum_{k=n}^{n_0-1}\sum_{j=1}^{M_k}
	L_{\TT_{k,j}}f=\sum_{k=n}^{n_0-1}\sum_{j=1}^{M_k}
	L_{\TT_{k,j}}f,
\end{aligned}
\]
which gives the desired decomposition in (v). The proof is complete.
\end{proof}

\medskip

\section{An $L^p$ bound for the linearized operator}
\label{sec:tree-variation}\label{sec:summation}\label{sec:full}

In this section, we prove Theorem~\ref{thm:fixed-linearization}. For this purpose, we first derive a martingale variation inequality from $(\Wq)$
and use it to estimate a single tree. Combining the tree estimate
with the decomposition of Section~\ref{sec:trees}, we then establish
restricted estimates for fixed linearizations. Finally, we remove the
major-subset restriction and apply interpolation to obtain the
$L^p$ bound for the linearized operator.  Let us begin with the following martingale variation inequality.

\begin{lemma}\label{lem:mart-variation}
Let  $2\le q<\infty$.  Suppose that $X$ is a Banach space satisfying $(\Wq)$ with constant
$A_q$. Then  for every $1<a<\infty$
and $q<r<\infty$, we have
\begin{equation*}\label{eq:lepingle}
\norm{\Var^r(\E_kf:k\in\Z)}_{L^a(\R_+)}
\le C_{a,r,q,X}\norm f_{L^a(\R_+;X)},
\qquad f\in L^a(\R_+;X).
\end{equation*}
\end{lemma}

\begin{proof}
	We first show that $(\Wq)$ implies dyadic martingale
	cotype $q$. Fix an arbitrary finite dyadic martingale
	$(x_n)_{n=0}^N$ as defined in Section~3, where
	\[
	\Omega_N=\{-1,1\}^N,
	\qquad
	\mathscr F_n=\sigma(\omega_1,\ldots,\omega_n),
	\qquad
	x_n=\E(x_N\mid\mathscr F_n).
	\]
	For each $\omega\in\Omega_N$, let $I_\omega\subseteq[0,1)$
	be the dyadic interval of length $2^{-N}$ on which the
	first $N$ Rademacher functions take the values
	$\omega_1,\ldots,\omega_N$, respectively. Define an
	$X$-valued step function $f$ on $[0,1)$ by
	\[
	f(t)=x_N(\omega),\qquad t\in I_\omega.
	\]
	By the definition of dyadic conditional expectations,
	\[
	\E_nf(t)
	=\E(x_N\mid\mathscr F_n)(\omega)
	=x_n(\omega),
	\qquad t\in I_\omega,\quad 0\le n\le N.
	\]
	In particular,
	\[
	\norm f_{L^q([0,1);X)}^q
	=\E\norm{x_N}_X^q.
	\]
	For each dyadic interval $I\subseteq[0,1)$ with
	$|I|\ge2^{-N}$, set
	\[
	P(I)=I\times[2^{N+1}-2|I|^{-1},2^{N+1}),
	\qquad
	m=w_{2^{N+1}-1}.
	\]
	Then the down-tiles $P(I)_d$ are pairwise disjoint:
	their frequency intervals are disjoint across
	different scales, while their time intervals are
	disjoint at each fixed scale. By \eqref{eq:restriction},
	$
	mw_{P(I)_d}=\pm h_I.
	$
	Since $|m|=1$ on $[0,1)$, applying $(\Wq)$
	to $mf$, we get
	\begin{equation}\label{eq:dyadic-cotype}
		\sum_{\substack{I\subseteq[0,1)\\|I|\ge2^{-N}}}
		\norm{\ip f{h_I}h_I}_{L^q(X)}^q
		\le A_q^q\norm f_{L^q(X)}^q.
	\end{equation}
For each $1\le n\le N$, the Haar expansion of the
	dyadic martingale difference gives
	\[
	x_n-x_{n-1}
	=\E_nf-\E_{n-1}f
	=\sum_{\substack{I\subseteq[0,1)\\
			|I|=2^{-(n-1)}}}
	\ip f{h_I}h_I.
	\]
	Since the intervals at each fixed scale are
	pairwise disjoint, we have
	\[
	\E\norm{x_n-x_{n-1}}_X^q
	=
	\sum_{\substack{I\subseteq[0,1)\\
			|I|=2^{-(n-1)}}}
	\norm{\ip f{h_I}h_I}_{L^q(X)}^q.
	\]
	Consequently, it follows from \eqref{eq:dyadic-cotype} that
	\[
	\begin{aligned}
		\sum_{n=1}^N\E\norm{x_n-x_{n-1}}_X^q
		&=
		\sum_{\substack{I\subseteq[0,1)\\
				2^{-N}<|I|\le1}}
		\norm{\ip f{h_I}h_I}_{L^q(X)}^q\\
		&\le A_q^q\norm f_{L^q(X)}^q=A_q^q\E\norm{x_N}_X^q.
	\end{aligned}
	\]
	Here the Haar coefficients on intervals of length
	$2^{-N}$ vanish, since $f$ is constant on those
	intervals. This shows that $X$ satisfies dyadic martingale
	cotype $q$.
	
	By \cite[Theorem~4.6]{Pisier},
	dyadic martingale cotype $q$ implies martingale
	cotype $q$ for arbitrary filtrations.   Therefore, $X$ satisfies martingale
	cotype $q$.
	

We now apply the martingale jump inequality
\cite[Theorem~1.3]{MSZ}. For a finite sequence
$(z_k)_{k\in\Lambda}$ in $X$, let
\[
N_\lambda(z_k:k\in\Lambda)
:=\sup\bigl\{m:\ k_0<\cdots<k_m,\
\norm{z_{k_i}-z_{k_{i-1}}}_X\ge\lambda,\
1\le i\le m\bigr\},
\]
where the indices belong to $\Lambda$ and an empty
collection of jumps has count zero. Since $X$ has martingale cotype $q$, applying  the jump inequality, for every
$1<b<\infty$ and finite $\Lambda\subseteq\Z$,
\[
\sup_{\lambda>0}
\norm{\lambda N_\lambda(\E_kf:k\in\Lambda)^{1/q}}_{L^b}
\le C_{b,q,X}\norm f_{L^b(X)},\qquad f\in L_b(X).
\]
By \cite[Lemma~2.12]{MSZ}, this jump estimate yields,
for every $r>q$,
\[
\bigl|\{\Var^r(\E_kf:k\in\Lambda)>\lambda\}\bigr|
\le C_{b,r,q,X}
\lambda^{-b}\norm f_{L^b(X)}^b,
\qquad \lambda>0.
\]
Then, applying
the Marcinkiewicz interpolation theorem
\cite[Lemma~2.5]{GLY} at two exponents
$1<b_0<a<b_1<\infty$, we obtain
\[
\norm{\Var^r(\E_kf:k\in\Lambda)}_{L^a}
\le C_{a,r,q,X,A_q}\norm f_{L^a(X)},
\]
uniformly over finite $\Lambda\subseteq\Z$. Finally, take $\Lambda=\{-m,\ldots,m\}$ and let
$m\to\infty$. Since the corresponding variation
seminorms increase to the full seminorm, monotone
convergence proves the desired result.
\end{proof}

We next combine the martingale variation inequality of
Lemma~\ref{lem:mart-variation} with the size and density
estimates to bound the linearized operator associated
with a single tree.

\begin{lemma}\label{lem:tree}
	Suppose that $2\le q<r<\infty$ and that $X$ satisfies
	$(\Wq)$ with constant $A_q$. Let
	$\PP\subseteq\mathcal B$ be a finite bitile family
	and fix a linearization $L_{\PP}$.
	Let $\TT\subseteq\PP$ be a tree with top
	$T\in\mathcal B$. Then for every
	$f\in L^1_{\mathrm{loc}}(\R_+;X)$,
	$g\in L^{r'}(\R_+)$ and $1\le a\le r'$, we have
	\[
	\norm{gL_{\TT}f}_{L^a(\R_+)}
	\le C\size_q(\TT;f)\dens(\TT;g)|I_T|^{1/a},
	\]
	where $L_{\TT}$ is the restriction of
	$L_{\PP}$ to $\TT$.
\end{lemma}

\begin{proof}
	The estimate is immediate when $\TT=\varnothing$.
	Suppose that $\TT\ne\varnothing$. Since $P\le T$
	for every $P\in\TT$, the function $L_{\TT}f$ is
	supported on $I_T$. Hence, for $1\le a\le r'$,
	H\"older's inequality gives
	\[
	\norm{gL_{\TT}f}_{L^a}
	\le |I_T|^{1/a-1/r'}
	\norm{gL_{\TT}f}_{L^{r'}}.
	\]
	It therefore suffices to prove the estimate for
	$a=r'$. Write
	$\TT=\TT_u\mathbin{\dot\cup}\TT_d$. For convenience, we split the proof into several steps.
	
	\medskip
	\noindent\textit{Step 1. Localization of the density.}	Let $\mathcal J$ consist of the maximal dyadic
	subintervals of $I_T$ containing no $I_P$ with
	$P\in\TT$. For every $x\in I_T$, the dyadic
	subinterval of $I_T$ containing $x$ and having
	length
$\frac12\min_{P\in\TT}|I_P|$
	contains no such $I_P$. Taking its largest ancestor
	with the same property gives an interval in
	$\mathcal J$ containing $x$. By dyadic nesting,
	the members of $\mathcal J$ are pairwise disjoint.
	Moreover, each $J\in\mathcal J$ satisfies
$
	|J|\ge\frac12\min_{P\in\TT}|I_P|.
$
	Consequently, $\mathcal J$ is finite and
	\[
	I_T=\mathop{\dot\bigcup}_{J\in\mathcal J}J.
	\]
	Fix $J\in\mathcal J$. Since $\TT\ne\varnothing$,
	we have $J\subsetneq I_T$, and hence its dyadic
	parent $J^{(1)}$ is contained in $I_T$. By the
	maximality of $J$, there exists $Q(J)\in\TT$ such that $
	I_{Q(J)}\subseteq J^{(1)}.$
	Let $P(J)$ be the bitile with time interval $J^{(1)}$
	and frequency interval equal to the dyadic ancestor
	of $\omega_T$ of length $2/|J^{(1)}|$. Since
$$
	|\omega_T|\le\frac{2}{|J^{(1)}|}
	\le|\omega_{Q(J)}|,
$$
	and $\omega_{Q(J)}$ contains $\omega_T$, we have
	\[
	\omega_T\subseteq\omega_{P(J)}
	\subseteq\omega_{Q(J)}.
	\]
	Thus
	\[
	Q(J)\le P(J)\le T.
	\]
	For any $P\in\TT$ with $I_P\cap J\ne\varnothing$,
	the definition of $\mathcal J$ and dyadic nesting
ensure that
$
	J^{(1)}\subseteq I_P.
$
	Since both $\omega_P$ and $\omega_{P(J)}$ contain
	$\omega_T$, and
	\[
	|\omega_P|=\frac{2}{|I_P|}
	\le\frac{2}{|J^{(1)}|}
	=|\omega_{P(J)}|,
	\]
	we obtain
	\[
	\omega_P\subseteq\omega_{P(J)}.
	\]
	This common frequency interval allows us to apply
	the density estimate on every $J\in\mathcal J$.
	Indeed, taking $Q(J)$ and $P(J)$ in the two suprema
	defining \eqref{eq:density}, respectively, we get
	\begin{align}\label{eq:tree-local-density}
		&\int_J |g|^{r'}
		\sum_{j:N_j\in\omega_{P(J)}}
		\norm{a_j^*}_{X^*}^{r'}\dd x\notag\\
		&\qquad\le
		\int_{I_{P(J)}} |g|^{r'}
		\sum_{j:N_j\in\omega_{P(J)}}
		\norm{a_j^*}_{X^*}^{r'}\dd x\\
		&\qquad\le
		|I_{P(J)}|\dens(\TT;g)^{r'}
		=2|J|\dens(\TT;g)^{r'}.\notag
	\end{align}
	
	\medskip
\noindent\textit{Step 2. Estimate for the down-tree $\TT_d$.} 	Fix $x\in\R_+$. There is at most one pair
	$(P,j)$ with $P\in\TT_d$, $1\le j\le K$, and
	$x\in I_P\cap A(P,j)$. Indeed, a fixed bitile can
	contribute to at most one index $j$, by the proof
	of Proposition~\ref{prop:exact-linearization}.
	Suppose that two distinct pairs $(P,j)$ and
	$(P',j')$ contribute. Then $P\ne P'$. Since
	$\TT_d$ is a down-tree,
	\[
	\varnothing\ne\{x\}\times\omega_{T_d}
	\subseteq P\cap P'.
	\]
	By Subsection~\ref{subsec:tile-orders}, the two
	bitiles are comparable. Since they are distinct,
	their time intervals have different lengths.
	After interchanging the pairs if necessary, assume
	that $|I_P|<|I_{P'}|$. As the time intervals are
	dyadic, we have
	\[
	|\omega_{P'}|=\frac{2}{|I_{P'}|}
	\le\frac{1}{|I_P|}
	=|\omega_{P_d}|.
	\]
	Both frequency intervals contain $\omega_{T_d}$,
	so dyadic nesting gives
	\[
	\omega_{P'}\subseteq\omega_{P_d}.
	\]
	The defining conditions of $A(P,j)$ and
	$A(P',j')$ therefore imply
	\[
	N_{j'}(x)<\sup\omega_{P_d}
	=\inf\omega_{P_u}\le N_j(x),
	\]
	and hence $j'<j$. Consequently,
	\[
	\begin{aligned}
		\inf\omega_P
		&\le\inf\omega_{P'}<\inf\omega_{P'_u}
		\le N_{j'}(x)\le N_{j-1}(x)<N_j(x)<\sup\omega_P,
	\end{aligned}
	\]
	contrary to $N_{j-1}(x)\notin\Int\omega_P$.
	This proves the asserted uniqueness.
	
	Now fix $J\in\mathcal J$ and $x\in J$.
	If there is a contributing pair $(P,j)$, then
	$I_P$ meets $J$, and hence
$
	N_j(x)\in[\omega_{P_u})
	\subseteq\omega_P\subseteq\omega_{P(J)}.
$
	By \eqref{eq:linearization}, the uniqueness
	proved above and \eqref{eq:singleton}, we obtain
	\begin{align*}
		|g(x)L_{\TT_d}f(x)|
		&\le |g(x)|\,|I_P|^{-1/2}
		\norm{\ip f{w_{P_d}}}_X
		\norm{a_j^*(x)}_{X^*}\\
		&\le \size_q(\TT;f)|g(x)|
		\norm{a_j^*(x)}_{X^*}\\
		&\le \size_q(\TT;f)|g(x)|
		\left(
		\sum_{\ell:N_\ell(x)\in\omega_{P(J)}}
		\norm{a_\ell^*(x)}_{X^*}^{r'}
		\right)^{1/r'}.
	\end{align*}
	If there is no contributing pair, then
	$L_{\TT_d}f(x)=0$, so the last estimate holds
	throughout $J$.	Since $L_{\TT_d}f$ is supported on $I_T$, summing
	over $\mathcal J$ and using
	\eqref{eq:tree-local-density}, we obtain
	\begin{equation}\label{eq:tree-down-bound}
	\begin{split}
		\norm{gL_{\TT_d}f}_{L^{r'}}&=\left(\sum_{J\in\mathcal J}
		\int_J|gL_{\TT_d}f|^{r'}\dd x\right)^{1/r'}\\
		&\le\size_q(\TT;f)
		\left(
		\sum_{J\in\mathcal J}
		\int_J |g|^{r'}
		\sum_{\ell:N_\ell\in\omega_{P(J)}}
		\norm{a_\ell^*}_{X^*}^{r'}\dd x
		\right)^{1/r'}\\
		&\le\size_q(\TT;f)
		\left(
		2\dens(\TT;g)^{r'}
		\sum_{J\in\mathcal J}|J|
		\right)^{1/r'}\\
		&=2^{1/r'}\size_q(\TT;f)
		\dens(\TT;g)|I_T|^{1/r'}.
\end{split}
	\end{equation}
	
	\medskip
	\noindent\textit{Step 3. Estimate for the up-tree $\TT_u$.}
To deal with the up-tree $\TT_u$, we first use a common
modulation to turn the Walsh packet functions into
Haar functions, and show that for each $j$, the
corresponding packet sum can be written as a
difference of two dyadic conditional expectations.
	Take $Q_T=T_u$ in \eqref{eq:two-components} and let
	$m_T$ be the corresponding modulation. Set
	\[
	\begin{aligned}
		\phi
		&=m_T\sum_{P\in\TT_u}
		\epsilon_P\ip f{w_{P_d}}w_{P_d}=\sum_{P\in\TT_u}
		\epsilon_P\eta_P\ip f{w_{P_d}}h_{I_P},
	\end{aligned}
	\]
	where $\eta_P$ is the sign in
	\eqref{eq:two-components}.
		Fix $J\in\mathcal J$ and $x\in J$. For $k\in\Z$,
	let $P_k(x)$ be the bitile $P\in\TT_u$ with
	$x\in I_P$ and $|I_P|=2^{-k}$, if it exists.
	There is at most one such bitile for each $k$:
	any two of them would intersect on
	$\{x\}\times\omega_T$, and bitiles with equal
	time intervals that intersect must coincide.
	For $1\le j\le K$, define
	\[
	\Lambda(x)=\{k\in\Z:P_k(x)\text{ exists}\},
	\qquad
	\Lambda_j(x)=
	\{k\in\Lambda(x):x\in A(P_k(x),j)\}.
	\]
	For $P_k(x)\in\TT_u$, the frequency intervals
	$\omega_{P_k(x)}$ and $\omega_{(P_k(x))_u}$
	contain $\omega_T$ and $\omega_{T_u}$, respectively.
	Their lengths increase with $k$, so both families
	are nested. In particular, if
	$k_-\le k\le k_+$ belong to $\Lambda(x)$ and
	$k_-,k_+\in\Lambda_j(x)$, then
	\[
	\begin{aligned}
		&N_j(x)
		\in[\omega_{(P_{k_-}(x))_u})
		\subseteq[\omega_{(P_k(x))_u}),\\
		&N_{j-1}(x)
		\in\R_+\setminus\Int\omega_{P_{k_+}(x)}
		\subseteq\R_+\setminus\Int\omega_{P_k(x)}.
	\end{aligned}
	\]
	Thus $k\in\Lambda_j(x)$. Hence, whenever
	$\Lambda_j(x)\ne\varnothing$, we may set
	\[
	a_j(x)=\min\Lambda_j(x),
	\qquad
	b_j(x)=1+\max\Lambda_j(x),
	\]
	and obtain
	\[
	\Lambda_j(x)=\Lambda(x)\cap[a_j(x),b_j(x)).
	\]
	By the Haar expansion of dyadic conditional
	expectation differences, we have
	\[
	(\E_{k+1}-\E_k)\phi
	=
	\sum_{\substack{P\in\TT_u\\|I_P|=2^{-k}}}
	\epsilon_P\eta_P\ip f{w_{P_d}}h_{I_P},
	\qquad k\in\Z.
	\]
	At $x$, this sum vanishes whenever
	$k\notin\Lambda(x)$, since every term is supported
	on its time interval. If $\Lambda_j(x)=\varnothing$,
	there is no contributing bitile for this index $j$,
	and the corresponding packet sum vanishes. Thus,
	for every nonempty $\Lambda_j(x)$, we obtain
	\begin{equation}\label{eq:up-tree-selected-block}
		\begin{split}
			&m_T(x)\sum_{P\in\TT_u}
			\epsilon_P\ip f{w_{P_d}}w_{P_d}(x)
			\one_{A(P,j)}(x)\\
			&\qquad=
			\sum_{k\in\Lambda_j(x)}
			[(\E_{k+1}-\E_k)\phi](x)\\
			&\qquad=
			\sum_{k=a_j(x)}^{b_j(x)-1}
			[(\E_{k+1}-\E_k)\phi](x)\\
			&\qquad=
			(\E_{b_j(x)}\phi-\E_{a_j(x)}\phi)(x).
		\end{split}
	\end{equation}

	We now show that the intervals
	$[a_j(x),b_j(x))$ corresponding to nonempty
	$\Lambda_j(x)$ are pairwise ordered. Suppose that
	$j<j'$ and there exist
	$k\in\Lambda_j(x)$ and $k'\in\Lambda_{j'}(x)$
	with $k\le k'$. By the nesting of the frequency
	intervals,
	\[
	N_j(x)\in[\omega_{(P_k(x))_u})
	\subseteq\Int\omega_{P_k(x)}
	\subseteq\Int\omega_{P_{k'}(x)}.
	\]
	Since $N_{j'}(x)\in[\omega_{(P_{k'}(x))_u})$,
	we have
	\[
	\inf\omega_{P_{k'}(x)}<N_j(x)
	\le N_{j'-1}(x)<N_{j'}(x)
	<\sup\omega_{P_{k'}(x)},
	\]
	contrary to $x\in A(P_{k'}(x),j')$. Therefore,
	\begin{equation}\label{eq:up-tree-level-order}
		\max\Lambda_{j'}(x)<\min\Lambda_j(x),
		\qquad j<j'
	\end{equation}
	whenever both sets are nonempty. In particular,
	\[
	a_{j'}(x)<b_{j'}(x)\le a_j(x)<b_j(x),
	\qquad j<j'.
	\]

We next show that the corresponding martingale
increments are controlled by the dyadic $r$-variation
of $\phi$ on each interval $J\in\mathcal J$.
Note that every Haar function $h_{I_P}$ in $\phi$ is constant
on $J$: if $I_P$ meets $J$, then
$J^{(1)}\subseteq I_P$, so $J$ lies in one of the
two children of $I_P$; otherwise $h_{I_P}$ vanishes
on $J$. Moreover, by the definition of dyadic
conditional expectations,
\[
\E_k h_{I_P}=
\begin{cases}
	0,&2^{-k}\ge |I_P|,\\
	h_{I_P},&2^{-k}<|I_P|.
\end{cases}
\]
Consequently, $\E_k\phi$ is constant on $J$ for
every $k\in\Z$. Since $\phi$ is a finite Haar sum,
$\E_k\phi$ changes only finitely many times.
We may therefore define
\[
V_J=\Var^r(\E_k\phi(x):k\in\Z),
\qquad x\in J,
\]
which is independent of the choice of $x\in J$. For each fixed $x\in J$, if there are nonempty
intervals $[a_j(x),b_j(x))$, list all their endpoints
in increasing order as
$\ell_0<\cdots<\ell_M$, counting each shared endpoint
only once. If there are no such intervals, set $M=0$
and interpret the sums below as zero. By
\eqref{eq:up-tree-level-order}, the endpoints of each
nonempty interval are consecutive in this list. Hence
\begin{equation}\label{eq:up-tree-block-variation}
	\begin{split}
		&\left(
		\sum_{j:\Lambda_j(x)\ne\varnothing}
		\norm{(\E_{b_j(x)}\phi-\E_{a_j(x)}\phi)(x)}_X^r
		\right)^{1/r}\\
		&\qquad\le
		\left(
		\sum_{i=1}^M
		\norm{(\E_{\ell_i}\phi-\E_{\ell_{i-1}}\phi)(x)}_X^r
		\right)^{1/r}
		\le V_J.
	\end{split}
\end{equation}

For each $j$ with $\Lambda_j(x)\ne\varnothing$,
by the definition of $\Lambda_j(x)$ and the nesting
of the frequency intervals, we have
\[
N_j(x)\in[\omega_{(P_k(x))_u})
\subseteq\omega_{P_k(x)}
\subseteq\omega_{P(J)},
\qquad k\in\Lambda_j(x).
\]
In particular,
\begin{equation}\label{eq:inclusion}
\{j:\Lambda_j(x)\ne\varnothing\}
\subseteq
\{j:N_j(x)\in\omega_{P(J)}\}.
\end{equation}
Using $|m_T(x)|=1$ on $J$, the definition of linearization
operators and \eqref{eq:up-tree-selected-block}, we get
\begin{align*}
	|g(x)L_{\TT_u}f(x)|
	&=|g(x)|\,|m_T(x)L_{\TT_u}f(x)|\\
	&=|g(x)|
	\left|
	\sum_{j:\Lambda_j(x)\ne\varnothing}
	\ip{
		m_T(x)\sum_{P\in\TT_u}
		\epsilon_P\ip f{w_{P_d}}w_{P_d}(x)
		\one_{A(P,j)}(x)
	}
	{a_j^*(x)}
	\right|\\
	&=|g(x)|
	\left|
	\sum_{j:\Lambda_j(x)\ne\varnothing}
	\ip{(\E_{b_j(x)}\phi-\E_{a_j(x)}\phi)(x)}
	{a_j^*(x)}
	\right|.
\end{align*}
Applying H\"older's inequality,
\eqref{eq:up-tree-block-variation} and the
inclusion \eqref{eq:inclusion}, we arrive at
\begin{align*}
	|g(x)L_{\TT_u}f(x)|
	&\le |g(x)|
	\left(
	\sum_{j:\Lambda_j(x)\ne\varnothing}
	\norm{(\E_{b_j(x)}\phi-\E_{a_j(x)}\phi)(x)}_X^r
	\right)^{1/r}\\
	&\qquad{}\times
	\left(
	\sum_{j:\Lambda_j(x)\ne\varnothing}
	\norm{a_j^*(x)}_{X^*}^{r'}
	\right)^{1/r'}\\
	&\le |g(x)|V_J
	\left(
	\sum_{j:N_j(x)\in\omega_{P(J)}}
	\norm{a_j^*(x)}_{X^*}^{r'}
	\right)^{1/r'},
	\qquad x\in J.
\end{align*}
Integrating the pointwise estimate over the partition
$I_T=\mathop{\dot\bigcup}_{J\in\mathcal J}{J}$ and using the fact that
$V_J$ is constant on each $J$, we obtain, by
\eqref{eq:tree-local-density},
\begin{align*}
	\norm{gL_{\TT_u}f}_{L^{r'}}
	&=
	\left(
	\sum_{J\in\mathcal J}
	\int_J|gL_{\TT_u}f|^{r'}\dd x
	\right)^{1/r'}\\
	&\le
	\left(
	\sum_{J\in\mathcal J}V_J^{r'}
	\int_J|g|^{r'}
	\sum_{j:N_j\in\omega_{P(J)}}
	\norm{a_j^*}_{X^*}^{r'}\dd x
	\right)^{1/r'}\\
	&\le
	2^{1/r'}\dens(\TT;g)
	\left(
	\sum_{J\in\mathcal J}|J|V_J^{r'}
	\right)^{1/r'}\\
	&=
	2^{1/r'}\dens(\TT;g)
	\norm{\Var^r(\E_k\phi:k\in\Z)}_{L^{r'}(I_T)}\\
	&\le
	C\dens(\TT;g)\norm{\phi}_{L^{r'}(\R_+;X)},
\end{align*}
where the last inequality follows from
Lemma~\ref{lem:mart-variation}. It remains to estimate the $L^{r'}$ norm of $\phi$.
The condition $r>q\ge2$ ensures that $r'<q$.
Moreover, $\phi$ is supported on $I_T$ and
$|m_T|=1$  on $I_T$. Thus, by H\"older's inequality
and a decomposition according to the signs
$\epsilon_P$, we obtain
\begin{equation}
	\begin{split}
		\norm{\phi}_{L^{r'}(\R_+;X)}
		&\le
		|I_T|^{1/r'-1/q}
		\left\|
		\sum_{P\in\TT_u}
		\epsilon_P\ip f{w_{P_d}}w_{P_d}
		\right\|_{L^q(\R_+;X)}\\
		&\le
		|I_T|^{1/r'-1/q}
		\sum_{\sigma\in\{-1,1\}}
		\left\|
		\sum_{\substack{P\in\TT_u\\\epsilon_P=\sigma}}
		\ip f{w_{P_d}}w_{P_d}
		\right\|_{L^q(\R_+;X)}\\
		&\le
		2\size_q(\TT;f)|I_T|^{1/r'}.
		\label{eq:signed-tree-size}
	\end{split}
\end{equation}
Indeed, each nonempty sign class is a subfamily
of the up-tree $\TT_u$ with top $T$, so the last
inequality follows directly from the definition
of size in \eqref{eq:size}.
 Combining
	\eqref{eq:signed-tree-size} with the preceding estimate,
	we obtain
	\begin{equation}\label{eq:tree-up-bound}
		\norm{gL_{\TT_u}f}_{L^{r'}}
		\le
		C\size_q(\TT;f)\dens(\TT;g)|I_T|^{1/r'}.
	\end{equation}
	
	\medskip
	\noindent\textit{Step 4. Completion of the proof.}	
 Combining
\eqref{eq:tree-down-bound} with \eqref{eq:tree-up-bound},
we conclude
\[
\begin{aligned}
	\norm{gL_{\TT}f}_{L^{r'}}
	&\le
	\norm{gL_{\TT_d}f}_{L^{r'}}
	+\norm{gL_{\TT_u}f}_{L^{r'}}\\
	&\le C\size_q(\TT;f)\dens(\TT;g)|I_T|^{1/r'},
\end{aligned}
\]
where $C$ depends only on $q,r,X$ and $A_q$.
This proves the estimate for $a=r'$, and hence
for the full range $1\le a\le r'$.
\end{proof}

%

We also need the following lemma which converts a packing condition on
dyadic intervals into an $L^b$ bound for their
overlap.

\begin{lemma}\label{lem:counting}
Let $(I_j)_{j=1}^N$ be a finite family of dyadic
intervals (allowing repetitions). Suppose that for some
$D>0$,
\begin{equation}\label{eq:hypo-pack}
\sum_{\substack{1\le j\le N\\I_j\subseteq J}}|I_j|
\le D|J|,
\qquad J\in\D.
\end{equation}
Then for every $1\le b<\infty$,
\[
\left\|\sum_{j=1}^N\one_{I_j}\right\|_{L^b}
\le C_b\left(\sum_{j=1}^N|I_j|\right)^{1/b}
D^{1-1/b}.
\]
\end{lemma}

\begin{proof}
The result is immediate if the family $(I_j)_{j=1}^N$ is empty.
Now set $f=\sum_{j=1}^N\one_{I_j}$. Fix $1\le b<\infty$ and choose an integer
$m>b$ with $m\ge2$. For any $m$ dyadic intervals with nonempty
intersection, their intersection equals a
shortest one, which is contained in every
other interval. Summing over the $m$ possible
choices of this shortest interval and applying
the packing hypothesis \eqref{eq:hypo-pack}, we obtain
\[
\begin{aligned}
\int f^m
&=\sum_{j_1,\ldots,j_m=1}^N
|I_{j_1}\cap\cdots\cap I_{j_m}|\\
&\le m\sum_{j_2,\ldots,j_m=1}^N
\sum_{\substack{1\le j_1\le N\\
I_{j_1}\subseteq\bigcap_{\nu=2}^m I_{j_\nu}}}
|I_{j_1}|\\
&\le mD\sum_{j_2,\ldots,j_m=1}^N
\left|\bigcap_{\nu=2}^m I_{j_\nu}\right|\\
&=mD\int f^{m-1}.
\end{aligned}
\]
Iterating this estimate and using
$\int f=\sum_{j=1}^N|I_j|$, we obtain
\[
\int f^m \le m!D^{m-1}\int  k  = m!D^{m-1} \sum_{j=1}^N|I_j|.
\]
By H\"older's inequality,
\[
\begin{aligned}
\norm{f}_{L^b}^b
&\le
\left(\int f\right)^{(m-b)/(m-1)}
\left(\int f^m\right)^{(b-1)/(m-1)}\\
&\le
(m!)^{(b-1)/(m-1)}
D^{b-1}\sum_{j=1}^N|I_j|.
\end{aligned}
\]
Taking $b$th roots completes the proof.
\end{proof}

The next result establishs a major-subset estimate for the
linearized operators, which will be used to prove
their $L^p$ boundedness.


%
%
%

\begin{proposition}\label{prop:major-subsets}
Let $2\le q<\infty$, $1<p<\infty$ and
$\max\{q,p'(q-1)\}<r<\infty$. Suppose that $X$
is a complex UMD space satisfying $(\Wq)$
with constant $A_q$.
For any measurable sets $F$, $E\subseteq\R_+$
with $0<|F|,|E|<\infty$, there exist measurable
subsets $F'\subseteq F$ and $E'\subseteq E$ such that
\[
\begin{cases}
F'=F,\quad |E'|\ge |E|/2,
& \text{if } |F|\le |E|,\\
E'=E,\quad |F'|\ge |F|/2,
& \text{if } |F|>|E|.
\end{cases}
\]
Moreover, for every finite model linearization $L_{\PP}$,
every strongly measurable $f:\R_+\to X$ and
measurable $g:\R_+\to\mathbb C$ satisfying
$
\norm{f}_X\le\one_{F'},$
$|g|\le\one_{E'}
$
almost everywhere, we have
\[
\left|\int_{\R_+}g(x)L_{\PP}f(x)\dd x\right|
\le C_{p,r,q,X,A_q}|F|^{1/p}|E|^{1/p'}.
\]
\end{proposition}

\begin{proof}

\noindent\textit{Case 1: $|F|\le |E|$.}
Let $\Md$ denote the dyadic Hardy--Littlewood
maximal operator. Set
\[
G=\{\Md\one_F>4|F|/|E|\},
\qquad F'=F,\qquad E'=E\setminus G.
\]
By the weak-type $(1,1)$ inequality for $\Md$
(see e.g. \cite[Exercise~2.1.12(a)]{GrafakosCFA}),
\[
|G|\le (4|F|/|E|)^{-1}|F|
=\frac{|E|}{4}.\]
Hence,
$$|E'|\ge |E|-|G|\ge\frac{3|E|}{4}\geq \frac{|E|}{2}.
$$

Fix $L_{\PP}$ and
functions $f,g$ satisfying $\norm{f}_X\le\one_F$ and $|g|\le\one_{E'}$
almost everywhere. Set
\[
\PP_G=\{P\in\PP:I_P\nsubseteq G\}.
\]
For every $P\in\PP\setminus\PP_G$, the packet
$w_{P_d}$ is supported on $I_P\subseteq G$,
where $g=0$ almost everywhere. Thus,
\begin{equation}\label{eq:useful-identity}
	gL_{\PP}f=gL_{\PP_G}f.
\end{equation}
Fix $R\in(\PP_G)^{\cvx}$. By the definition of
the convex hull, there exist $P_0,Q_0\in\PP_G$
such that $P_0\le R\le Q_0$. Since
$I_{P_0}\nsubseteq G$, we may choose
$x_R\in I_{P_0}\setminus G\subseteq I_R$.
It follows that, for $\nu\in\{d,u\}$,
\[
\begin{aligned}
	|I_R|^{-1/2}\norm{\ip f{w_{R_\nu}}}_X
	&\le\frac1{|I_R|}
	\int_{I_R}\norm{f(x)}_X\dd x\le\Md\one_F(x_R)
	\le\frac{4|F|}{|E|}.
\end{aligned}
\]
Applying Lemma~\ref{prop:size-bound} with
$M=4|F|/|E|$, we obtain
\begin{equation}\label{eq:small-F-size}
	\size_q(\PP_G;f)
	\le8\beta_{q,X}\frac{|F|}{|E|}.
\end{equation}
Moreover, since
$|g|\le\one_E$, it is clear that
\begin{equation}\label{eq:small-F-dense}
\dens(\PP_G;g)\le1.
\end{equation}
Now apply Theorem~\ref{prop:joint-levels} to
$\PP_G$.
By \eqref{eq:useful-identity} and
Theorem~\ref{prop:joint-levels}\textup{(v)}, we have
\[
\int_{\R_+}gL_{\PP}f\dd x
=\int_{\R_+}gL_{\PP_G}f\dd x
=\sum_{n,j}\int_{\R_+}gL_{\TT_{n,j}}f\dd x.
\]
Using \eqref{eq:level-size}, \eqref{eq:level-density},
\eqref{eq:small-F-size} and \eqref{eq:small-F-dense}, we get
\[
\begin{aligned}
	\size_q(\TT_{n,j};f)
	&\le C\min\left\{
	\frac{|F|}{|E|},(2^n|F|)^{1/q}
	\right\},\\
	\dens(\TT_{n,j};g)
	&\le C\min\left\{
	1,(2^n|E|)^{1/r'}
	\right\}\le C(2^n|E|)^{1/r'}.
\end{aligned}
\]
Hence, applying Lemma~\ref{lem:tree} with $a=1$,
the total packing estimate \eqref{eq:level-total} and the above estimates, we obtain
\begin{align}	\label{eq:small-F-sum}
	\left|\int_{\R_+}gL_{\PP}f\dd x\right|
	&\le\sum_{n,j}\norm{gL_{\TT_{n,j}}f}_{L^1}
	\notag\\
	&\le C\sum_{n,j}
	\size_q(\TT_{n,j};f)
	\dens(\TT_{n,j};g)|I_{T_{n,j}}|\\
	&\le C|E|^{1/r'}
	\sum_{n\in\Z}2^{-n/r}
	\min\left\{
	\frac{|F|}{|E|},(2^n|F|)^{1/q}
	\right\}.\notag
\end{align}
Split the sum at
$2^n=|F|^{q-1}/|E|^q$. Since $r>q$, both
geometric series converge and
\[
\begin{aligned}
	&|E|^{1/r'}
	\sum_{n\in\Z}2^{-n/r}
	\min\left\{
	\frac{|F|}{|E|},(2^n|F|)^{1/q}
	\right\}\\
	&\quad\le
	|F|^{1/q}|E|^{1/r'}
	\sum_{\substack{n\in\Z\\
			2^n\le |F|^{q-1}/|E|^q}}
	2^{n(1/q-1/r)}\\
	&\qquad+
	|F||E|^{-1/r}
	\sum_{\substack{n\in\Z\\
			2^n>|F|^{q-1}/|E|^q}}
	2^{-n/r}\\
	&\quad\le C_{q,r}
	|F|^{1-(q-1)/r}|E|^{(q-1)/r}.
\end{aligned}
\]
Combining this with \eqref{eq:small-F-sum}, we get
\[
\begin{aligned}
	\left|\int_{\R_+}gL_{\PP}f\dd x\right|
	&\le C|F|^{1-(q-1)/r}|E|^{(q-1)/r}\\
	&=C|F|^{1/p}|E|^{1/p'}
	\left(\frac{|F|}{|E|}\right)^{
		1/p'-(q-1)/r}\\
	&\le C|F|^{1/p}|E|^{1/p'},
\end{aligned}
\]
where the last inequality follows from
$r>p'(q-1)$ and $|F|\le|E|$.

\medskip

\medskip
\noindent\textit{Case 2. $|F|>|E|$.} Put
\[
G=\{\Md\one_E>8|E|/|F|\},
\qquad F'=F\setminus G,\qquad E'=E.
\]
By the weak-type $(1,1)$ inequality for $\Md$,
\[
|G|\le\frac{|F|}{8},
\qquad
|F'|\ge|F|-|G|\ge\frac{7|F|}{8}
\ge\frac{|F|}{2}.
\]

Fix $L_{\PP}$ and
functions $f,g$ satisfying
$\norm{f}_X\le\one_{F'}$ and $|g|\le\one_E$.
Set
\[
\PP_G=\{P\in\PP:I_P\nsubseteq G\}.
\]
For every $P\in\PP\setminus\PP_G$, we have
$I_P\subseteq G$, where $f=0$ almost everywhere.
Since $w_{P_d}$ is supported on $I_P$, it follows that
\[
\ip f{w_{P_d}}=0,
\qquad P\in\PP\setminus\PP_G.
\]
Thus $L_{\PP}f=L_{\PP_G}f$. Applying
Theorem~\ref{prop:joint-levels} to $\PP_G$
with the original majorants $F$, $E$, we obtain
\[
L_{\PP}f=L_{\PP_G}f
=\sum_{n,j}L_{\TT_{n,j}}f.
\]

We first improve the local packing estimate for the
top intervals. Fix a level $n$. If $J\in\D$ contains
a top interval $I_{T_{n,j}}$, choose
$P\in\TT_{n,j}\subseteq\PP_G$. Then
$I_P\subseteq I_{T_{n,j}}\subseteq J$ and
$I_P\nsubseteq G$, so there exists
$x\in I_P\setminus G\subseteq J\setminus G$.
It follows that
\[
\frac{|E\cap J|}{|J|}
\le\Md\one_E(x)
\le\frac{8|E|}{|F|},
\]
and hence
\[
\frac{|F\cap J|}{|F|}
+\frac{|E\cap J|}{|E|}
\le\frac{9|J|}{|F|}.
\]
Combining this with \eqref{eq:level-local}, we get
\begin{equation}\label{eq:improved-pack}
\sum_{j:I_{T_{n,j}}\subseteq J}|I_{T_{n,j}}|
\le C2^{-n}|F|^{-1}|J|,
\qquad J\in\D,
\end{equation}
which also holds trivially when $J$ contains no
top interval at level $n$. In particular, each
top interval satisfies
\[
|I_{T_{n,j}}|
\le
\sum_{i:I_{T_{n,i}}\subseteq I_{T_{n,j}}}
|I_{T_{n,i}}|
\le C2^{-n}|F|^{-1}|I_{T_{n,j}}|.
\]
Thus every nonempty level satisfies
\begin{equation}\label{eq:nonempty-level}
2^n|F|\le C.
\end{equation}

We now fix a nonempty level $n$ and estimate
the sum over its trees. Set
\[
k_n=\sum_j\one_{I_{T_{n,j}}},
\qquad t=\max\{r,p\}.
\]
Note that $L_{\TT_{n,j}}f$ is supported on
$I_{T_{n,j}}$. By H\"older's inequality, we obtain
\begin{equation}\label{eq:estimate-level-n}
\begin{split}
\left|\sum_j\int_{\R_+}
gL_{\TT_{n,j}}f\dd x\right|&\le
\int_{\R_+}k_n^{1/r}
\left(\sum_j|gL_{\TT_{n,j}}f|^{r'}\right)^{1/r'}
\dd x\\
&\le
\norm{k_n^{1/r}}_{L^t}
\left\|
\left(\sum_j|gL_{\TT_{n,j}}f|^{r'}\right)^{1/r'}
\right\|_{L^{t'}}.
\end{split}
\end{equation}
To estimate the first factor,  note that $t/r\ge1$. It follows from Lemma~\ref{lem:counting},
\eqref{eq:improved-pack} and
\eqref{eq:level-total} that
\begin{align}\label{eq:kn}
\norm{k_n^{1/r}}_{L^t}
&=\norm{k_n}_{L^{t/r}}^{1/r}\notag\\
&\le C
\left(\sum_j|I_{T_{n,j}}|\right)^{1/t}
\left(2^{-n}|F|^{-1}\right)^{1/r-1/t}\\
&\le C2^{-n/r}|F|^{1/t-1/r}.\notag
\end{align}
For the second factor, note that $
\left(\sum_j|gL_{\TT_{n,j}}f|^{r'}\right)^{1/r'}$
is supported on $E$. Since $t'\le r'$, we deduce from
H\"older's inequality, Lemma~\ref{lem:tree}
with $a=r'$,
\eqref{eq:level-size}, \eqref{eq:level-density} and \eqref{eq:level-total} that
\begin{equation}\label{eq:second-factor}
\begin{split}
&\left\|
\left(\sum_j|gL_{\TT_{n,j}}f|^{r'}\right)^{1/r'}
\right\|_{L^{t'}}\\
&\quad\le
|E|^{1/r-1/t}
\left(\sum_j
\norm{gL_{\TT_{n,j}}f}_{L^{r'}}^{r'}
\right)^{1/r'}\\
&\quad\le
C|E|^{1/r-1/t}
(2^n|F|)^{1/q}(2^n|E|)^{1/r'}
\left(\sum_j|I_{T_{n,j}}|\right)^{1/r'}\\
&\quad\le
C(2^n|F|)^{1/q}|E|^{1-1/t}.
\end{split}
\end{equation}
Plugging the estimates \eqref{eq:kn}, \eqref{eq:second-factor} into \eqref{eq:estimate-level-n},
we obtain
\begin{equation*}\label{eq:large-level}
\left|\sum_j\int_{\R_+}
gL_{\TT_{n,j}}f\dd x\right|
\le C2^{n(1/q-1/r)}
|F|^{1/t-1/r+1/q}|E|^{1-1/t}.
\end{equation*}

Finally, summing the above estimate
over the levels permitted by
\eqref{eq:nonempty-level} and using $r>q$,
we conclude
\[
\begin{aligned}
\left|\int_{\R_+}gL_{\PP}f\dd x\right|
&\le C|F|^{1/t-1/r+1/q}|E|^{1-1/t}
\sum_{\substack{n\in\Z\\2^n\le C/|F|}}
2^{n(1/q-1/r)}\\
&\le C|F|^{1/t}|E|^{1-1/t}\\
&=C|F|^{1/p}|E|^{1/p'}
\left(\frac{|F|}{|E|}\right)^{1/t-1/p}\\
&\le C|F|^{1/p}|E|^{1/p'},
\end{aligned}
\]
where the last inequality follows from
$t\ge p$ and $|F|>|E|$.
\end{proof}

We are now ready to prove Theorem~\ref{thm:fixed-linearization}.

\begin{proof}[Proof of Theorem~\ref{thm:fixed-linearization}]
Fix the parameters in Theorem~\ref{thm:fixed-linearization}
and choose a constant $A_q$ in $(\Wq)$ depending
only on $q$ and $X$. We divide the proof into several steps.

\medskip
\noindent\textit{Step 1. Preliminary bound and restricted-type constants.} We first establish an elementary bound for $L_{\PP}$.
As in the proof of
Proposition~\ref{prop:exact-linearization},
for each $P\in\PP$ and $x\in\R_+$,
at most one index $j$ satisfies $x\in A(P,j)$.
Together with \eqref{eq:directions}, this implies
\[
|L_{\PP}f(x)|
\le\sum_{P\in\PP}
\norm{\ip f{w_{P_d}}}_X|w_{P_d}(x)|.
\]
Since $|w_{P_d}|=|I_P|^{-1/2}\one_{I_P}$,
we have
\[
\norm{w_{P_d}}_{L^s}
\norm{w_{P_d}}_{L^{s'}}=1,
\qquad 1<s<\infty.
\]
Hence, by H\"older's inequality,
\begin{align}\label{eq:finite-linearization-bound}
\norm{L_{\PP}f}_{L^s}
&\le\sum_{P\in\PP}
\norm{\ip f{w_{P_d}}}_X
\norm{w_{P_d}}_{L^s}
\notag\\
&\le\sum_{P\in\PP}
\norm{f\one_{I_P}}_{L^s(X)}
\norm{w_{P_d}}_{L^{s'}}
\norm{w_{P_d}}_{L^s}\\
&\le|\PP|\norm{f}_{L^s(X)}.\notag
\end{align}
Since $r>\max\{q,p'(q-1)\}$, we may choose
\[
1<p_0<p<p_1<\infty,
\qquad
r>\max\{q,p_i'(q-1)\},
\quad i=0,1.
\]
Let $U$ consist of
bounded measurable scalar functions on $\R_+$
that vanish almost everywhere outside sets of
finite measure. Fix a strongly measurable
$v:\R_+\to X$ with $\norm{v}_X\le1$
almost everywhere. The pairing
\[
(h,g)\longmapsto
\int_{\R_+}g(x)L_{\PP}(vh)(x)\dd x,
\qquad h,g\in U,
\]
defines a complex bilinear form. By
\eqref{eq:finite-linearization-bound}
and H\"older's inequality, the above integral
is absolutely convergent. For $i=0,1$, let $C_i$ be the least constant
such that
\[
\left|\int_{\R_+}gL_{\PP}(vh)\dd x\right|
\le C_i|F|^{1/p_i}|E|^{1/p_i'}
\]
for all measurable $F,E\subseteq\R_+$ with
$0<|F|,|E|<\infty$ and all $h,g\in U$
satisfying $|h|\le\one_F$ and $|g|\le\one_E$.
By \eqref{eq:finite-linearization-bound}
and H\"older's inequality,
\[
\begin{aligned}
\left|\int_{\R_+}gL_{\PP}(vh)\dd x\right|
&\le|\PP|\norm{h}_{L^{p_i}}
\norm{g}_{L^{p_i'}}\le|\PP||F|^{1/p_i}|E|^{1/p_i'}.
\end{aligned}
\]
Hence $C_i\le|\PP|<\infty$.

\medskip
\noindent\textit{Step 2. Removal of the major-subset restriction.} Fix $i\in\{0,1\}$ and $F,E,h,g$ as in Step 1.
Choose the major subsets $F', E'$ from
Proposition~\ref{prop:major-subsets}
with exponent $p_i$. If $|F|\le|E|$, then $F'=F$ and
$|E\setminus E'|\le|E|/2$. Splitting $g$
according to $E'$ and $E\setminus E'$
and applying Proposition~\ref{prop:major-subsets}, we obtain
\[
\begin{aligned}
\left|\int_{\R_+}gL_{\PP}(vh)\dd x\right|
&\le
\left|\int_{\R_+}g\one_{E'}
L_{\PP}(vh)\dd x\right|+
\left|\int_{\R_+}g\one_{E\setminus E'}
L_{\PP}(vh)\dd x\right|\\
&\le
C_{p_i,r,q,X,A_q}
|F|^{1/p_i}|E|^{1/p_i'}+
C_i|F|^{1/p_i}
|E\setminus E'|^{1/p_i'}\\
&\le
\bigl(C_{p_i,r,q,X,A_q}
+2^{-1/p_i'}C_i\bigr)
|F|^{1/p_i}|E|^{1/p_i'}.
\end{aligned}
\]
If $|F|>|E|$, then $E'=E$ and
$|F\setminus F'|\le|F|/2$. By the linearity
of $L_{\PP}$, we similarly have
\[
\begin{aligned}
\left|\int_{\R_+}gL_{\PP}(vh)\dd x\right|
&\le
\left|\int_{\R_+}gL_{\PP}
(vh\one_{F'})\dd x\right|+
\left|\int_{\R_+}gL_{\PP}
(vh\one_{F\setminus F'})\dd x\right|\\
&\le
C_{p_i,r,q,X,A_q}
|F|^{1/p_i}|E|^{1/p_i'}+
C_i|F\setminus F'|^{1/p_i}
|E|^{1/p_i'}\\
&\le
\bigl(C_{p_i,r,q,X,A_q}
+2^{-1/p_i}C_i\bigr)
|F|^{1/p_i}|E|^{1/p_i'}.
\end{aligned}
\]
Taking the supremum over all admissible
$F,E,h,g$, we obtain
\[
C_i\le C_{p_i,r,q,X,A_q}
+\max\{2^{-1/p_i},2^{-1/p_i'}\}C_i.
\]
Since $1<p_i<\infty$, we have
$
\max\{2^{-1/p_i},2^{-1/p_i'}\}<1.$
Recalling that $C_i<\infty$, we conclude
\[
C_i\le
\frac{C_{p_i,r,q,X,A_q}}
{1-\max\{2^{-1/p_i},2^{-1/p_i'}\}},
\qquad i=0,1.
\]

\medskip
\noindent\textit{Step 3. The $L^p$ bound on the class $U$.}
To apply the bilinear interpolation theorem on $\R$,
we define a bilinear form for functions $H,G$
on $\R$ by
\[
(H,G)\longmapsto
\int_{\R_+}
(G|_{\R_+})L_{\PP}(v(H|_{\R_+}))\dd x.
\]
Let $H, G:\R\to\mathbb C$ be bounded measurable
functions vanishing almost everywhere outside
sets of finite measure and suppose that
$|H|\le\one_F$ and $|G|\le\one_E$
for measurable $F, E\subseteq\R$ of finite measure.
The bounds from Step 2 give
\[
\begin{aligned}
\left|\int_{\R_+}
(G|_{\R_+})L_{\PP}(v(H|_{\R_+}))\dd x\right|& \le
C_i|F\cap\R_+|^{1/p_i}
|E\cap\R_+|^{1/p_i'}\\
&\le
C_i|F|^{1/p_i}|E|^{1/p_i'},
\qquad \qquad\qquad i=0,1.
\end{aligned}
\]
Hence,
\[
\left|\int_{\R_+}
(G|_{\R_+})L_{\PP}(v(H|_{\R_+}))\dd x\right|
\le \min\left\{
C_0|F|^{1/p_0}|E|^{1/p_0'},
C_1|F|^{1/p_1}|E|^{1/p_1'}
\right\}.
\]
Choose
$
\theta=
\frac{1/p-1/p_1}{1/p_0-1/p_1}\in(0,1).
$ We further get
\[
\begin{aligned}
&\left|\int_{\R_+}
(G|_{\R_+})L_{\PP}(v(H|_{\R_+}))\dd x\right|\le
C_0^\theta C_1^{1-\theta}
|F|^{\theta/p_0+(1-\theta)/p_1}
|E|^{\theta/p_0'+(1-\theta)/p_1'}.
\end{aligned}
\]
Note that the exponents of $|F|$ and $|E|$
are $1/p$ and $1/p'$, respectively.
As $\theta$ varies in $(0,1)$, we obtain
restricted-type bounds in a relatively open
neighborhood of $(1/p,1/p')$ along the line
$\alpha+\beta=1$, with constants bounded by
$\max\{C_0,C_1\}$. Then, applying the bilinear case of
\cite[Theorem~3.2, pp.~16--17]{ThieleBook}, we derive
\[
\begin{aligned}
&\left|\int_{\R_+}
(G|_{\R_+})L_{\PP}(v(H|_{\R_+}))\dd x\right|\le
C_{p,r,q,X,A_q}
\norm{H}_{L^p(\R)}
\norm{G}_{L^{p'}(\R)}.
\end{aligned}
\]
For any $h,g\in U$, let $H$ and $G$ be their zero
extensions to $\R$, respectively. Namely,
\[
H(x)=
\begin{cases}
h(x),&x\ge0,\\
0,&x<0,
\end{cases}
\qquad
G(x)=
\begin{cases}
g(x),&x\ge0,\\
0,&x<0.
\end{cases}
\]
Then $H|_{\R_+}=h$, $G|_{\R_+}=g$, $
\norm{H}_{L^p(\R)}=\norm{h}_{L^p(\R_+)}$ $
\norm{G}_{L^{p'}(\R)}=\norm{g}_{L^{p'}(\R_+)}.$
Applying the preceding bilinear inequality to
$H$ and $G$, we obtain
\begin{equation}\label{eq:bilinear-strong}
\begin{split}
\left|\int_{\R_+}gL_{\PP}(vh)\dd x\right|
&\le C_{p,r,q,X,A_q}
\norm{H}_{L^p(\R)}
\norm{G}_{L^{p'}(\R)}\\
&=C_{p,r,q,X,A_q}
\norm{h}_{L^p(\R_+)}
\norm{g}_{L^{p'}(\R_+)}.
\end{split}
\end{equation}
Since $U$ is dense in $L^{p'}(\R_+)$
and closed under complex conjugation, we infer from duality and \eqref{eq:bilinear-strong} that
\[
\begin{aligned}
\norm{L_{\PP}(vh)}_{L^p}
&=
\sup_{\substack{g\in U\\
\norm{g}_{L^{p'}}\le1}}
\left|\int_{\R_+}gL_{\PP}(vh)\dd x\right|\le C_{p,r,q,X,A_q}\norm{h}_{L^p}.
\end{aligned}
\]

\medskip
\noindent\textit{Step 4. Extension to $L^p(\R_+;X)$.}
Let $f\in L^p(\R_+;X)$ and set
$h=\norm{f}_X$. Define
\[
v(x)=
\begin{cases}
f(x)/h(x),&h(x)>0,\\
0,&h(x)=0.
\end{cases}
\]
Then $v$ is strongly measurable,
$\norm{v}_X\le1$ and $f=vh$.
For $n\ge1$, put
\[
h_n=h\one_{[0,n)}\one_{\{h\le n\}}\in U.
\]
Since $0\le h_n\le h$ and $h_n\to h$
almost everywhere, dominated convergence gives
\[
\norm{h_n-h}_{L^p}\longrightarrow0.
\]
By \eqref{eq:finite-linearization-bound},
\[
\begin{aligned}
\norm{L_{\PP}(vh_n)-L_{\PP}f}_{L^p}
&=\norm{L_{\PP}(v(h_n-h))}_{L^p}\le|\PP|\norm{h_n-h}_{L^p}
\longrightarrow0.
\end{aligned}
\]
Applying the bound on $U$ (i.e., Step 3) and passing
to the limit, we conclude that
\[
\begin{aligned}
\norm{L_{\PP}f}_{L^p}
&=\lim_{n\to\infty}
\norm{L_{\PP}(vh_n)}_{L^p}\\
&\le C_{p,r,q,X,A_q}
\lim_{n\to\infty}\norm{h_n}_{L^p}\\
&=C_{p,r,q,X,A_q}\norm{f}_{L^p(X)}.
\end{aligned}
\]
The proof is complete.
\end{proof}

\medskip

\section{The UMD property and Walsh tile type}\label{sec:consequences}

This section is devoted to the proof of Theorem \ref{cor:tile-type}. The proof combines the known implication from finite Walsh tile type to the UMD property with a converse argument based on the disjoint-packet estimate and the finite-model bound established in Theorem~\ref{thm:model}.

\begin{proof}[Proof of Theorem \ref{cor:tile-type}]
If $X$ has Walsh tile type $\tau$ for some
$2\le\tau<\infty$, then it follows from \cite[Proposition~4.6]{HLP} that $X$ is UMD.

Conversely, suppose that $X$ is UMD.
By Theorem~\ref{prop:UMD-Wq}, there exists
$q_X\ge2$ such that $X$ satisfies
$(\mathrm W_{q_X})$ with some finite constant
$A_{q_X}$. Fix $q_X<\tau<\infty$,
$f\in L^\tau(\R_+;X)$ and a u-good family
$(\TT_j)_{j=1}^M$ of up-trees with corresponding
tops $T_j$. Set
$
\PP=\bigcup_{j=1}^M\TT_j
$
and take $\epsilon_P=1$ for every $P\in\PP$
in \eqref{eq:model}. For
$0\le j\le M$ and $x\in\R_+$, define
\begin{align*}
k(j,x)
&=\max\bigl(
\{k:1\le k\le j,\ x\in I_{T_k}\}\cup\{0\}
\bigr),\\
N_j(x)
&=
\begin{cases}
c(\omega_{T_{k(j,x)}}),&k(j,x)>0,\\
0,&k(j,x)=0.
\end{cases}
\end{align*}
For each $x$, the function $j\mapsto k(j,x)$
is nondecreasing. Hence
\eqref{eq:good-forest-centers} implies
\begin{equation}\label{eq:good-forest-monotone}
0=N_0(x)\le\cdots\le N_M(x).
\end{equation}
By \cite[Lemma~4.3(ii)]{HLP}, the defining
conditions of a u-good family give
\begin{equation*}
\begin{split}
\{P\in\TT_j:x\in I_P\}
=\bigl\{P\in\PP:\,
&x\in I_P,\ N_j(x)\in[\omega_{P_u}),N_{j-1}(x)\notin\Int\omega_P
\bigr\},
\end{split}
\end{equation*}
for $1\le j\le M$ and $x\in\R_+$.
Since $w_{P_d}$ is supported on $I_P$,
we obtain
\begin{equation}\label{eq:good-forest-packets}
\begin{split}
&\sum_{P\in\TT_j}
\ip f{w_{P_d}}w_{P_d}(x)=
\sum_{P\in\PP}
\ip f{w_{P_d}}w_{P_d}(x)
\one_{\left\{
\substack{
N_j(x)\in[\omega_{P_u})\\
N_{j-1}(x)\notin\Int\omega_P
}
\right\}}.
\end{split}
\end{equation}
Fix $x\in\R_+$. Since
$[\omega_{P_u})\subseteq\Int\omega_P$,
the sum in \eqref{eq:good-forest-packets}
vanishes whenever $N_{j-1}(x)=N_j(x)$.
Let $j_1<\cdots<j_K$ be the indices
$1\le j\le M$ for which
$N_{j-1}(x)<N_j(x)$, and put $j_0=0$.
If $K\ge1$, the sequence in
\eqref{eq:good-forest-monotone} is constant
between successive such indices. Thus
\begin{equation}\label{eq:good-forest-pair}
\bigl(N_{j_\ell-1}(x),N_{j_\ell}(x)\bigr)
=
\bigl(N_{j_{\ell-1}}(x),N_{j_\ell}(x)\bigr),
\qquad 1\le\ell\le K,
\end{equation}
and
\[
0=N_{j_0}(x)<N_{j_1}(x)<\cdots<N_{j_K}(x).
\]
Using these endpoints in \eqref{eq:model},
together with \eqref{eq:good-forest-packets}
and \eqref{eq:good-forest-pair}, we obtain
\[
\begin{aligned}
&\left(
\sum_{j=1}^M
\left\|
\sum_{P\in\TT_j}
\ip f{w_{P_d}}w_{P_d}(x)
\right\|_X^\tau
\right)^{1/\tau}=
\left(
\sum_{\ell=1}^K
\left\|
\sum_{P\in\TT_{j_\ell}}
\ip f{w_{P_d}}w_{P_d}(x)
\right\|_X^\tau
\right)^{1/\tau}\le C_{\tau,\PP}f(x).
\end{aligned}
\]
If $K=0$, all the tree sums vanish, so the
same inequality holds. Consequently,
\begin{equation*}
\sum_{j=1}^M
\left\|
\sum_{P\in\TT_j}
\ip f{w_{P_d}}w_{P_d}(x)
\right\|_X^\tau
\le \bigl(C_{\tau,\PP}f(x)\bigr)^\tau.
\end{equation*}
Since $\tau>q_X$, we have
$
\tau'(q_X-1)<\tau.
$
Thus the exponent conditions of
Theorem~\ref{thm:model} are satisfied with
$p=r=\tau$ and $q=q_X$.
Integrating the above estimate
and applying \eqref{eq:model-estimate},
we obtain
\[
\begin{aligned}
\sum_{j=1}^M
\left\|
\sum_{P\in\TT_j}
\ip f{w_{P_d}}w_{P_d}
\right\|_{L^\tau(\R_+;X)}^\tau& \le
\norm{C_{\tau,\PP}f}_{L^\tau(\R_+)}^\tau\\
&\le
\left(
\norm{C_{\tau,\PP}f}_{L^\tau(\R_+)}
+\norm{\widetilde C_{\tau,\PP}f}_{L^\tau(\R_+)}
\right)^\tau\\
&\le
C_\tau^\tau\norm{f}_{L^\tau(\R_+;X)}^\tau.
\end{aligned}
\]
This proves \eqref{eq:tile-type} and completes
the proof.
\end{proof}
\medskip

\section*{AI statement}
During the preparation of this work, the authors used GPT-6 Astra to assist with language editing and refinement of mathematical arguments. The author takes full responsibility for the content of this paper.



\begin{thebibliography}{99}

\bibitem{ADR}
A.~Amini-Harandi, I.~Doust, and G.~Robertson,
\emph{Roundness properties of Banach spaces},
J. Funct. Anal. \textbf{281} (2021), no.~10, article 109230.

\bibitem{AU}
A.~Amenta and G.~Uraltsev,
\emph{Banach-valued modulation invariant Carleson embeddings and outer-$L^p$ spaces: the Walsh case},
J. Fourier Anal. Appl. \textbf{26} (2020), no.~4, article 53, 54 pp.

\bibitem{AUVariation}
A.~Amenta and G.~Uraltsev,
\emph{Variational Carleson operators in UMD spaces},
J. Lond. Math. Soc. (2) \textbf{105} (2022), no.~3, 1363--1417.

\bibitem{Billard}
P.~Billard,
\emph{Sur la convergence presque partout des s\'eries de Fourier--Walsh des fonctions de l'espace $L^2(0,1)$},
Studia Math. \textbf{28} (1967), no.~3, 363--388.

\bibitem{Carleson}
L.~Carleson,
\emph{On convergence and growth of partial sums of Fourier series},
Acta Math. \textbf{116} (1966), 135--157.


\bibitem{CJWZ2026}
J.~Chang, Y.~Jiao, L.~Wu, and D.~Zhou,
\emph{A Carleson theorem in UMD spaces}, preprint, 2026.


\bibitem{DoLaceyWalsh}
Y.~Do and M.~T.~Lacey,
\emph{Weighted bounds for variational Walsh--Fourier series},
J. Fourier Anal. Appl. \textbf{18} (2012), no.~6, 1318--1339.

\bibitem{Fefferman}
C.~Fefferman,
\emph{Pointwise convergence of Fourier series},
Ann. of Math. (2) \textbf{98} (1973), no.~3, 551--571.

\bibitem{GrafakosCFA}
L.~Grafakos,
\emph{Classical Fourier Analysis},
3rd ed., Graduate Texts in Mathematics, vol.~249,
Springer, New York, 2014.

\bibitem{GLY}
L.~Grafakos, L.~Liu, and D.~Yang,
\emph{Vector-valued singular integrals and maximal functions on spaces
of homogeneous type},
Math. Scand. \textbf{104} (2009), no.~2, 296--310.

\bibitem{Hunt}
R.~A.~Hunt,
\emph{On the convergence of Fourier series},
Orthogonal Expansions and Their Continuous Analogues (Proc. Conf., Edwardsville, Ill., 1967),
Southern Illinois Univ. Press, Carbondale, IL, 1968, pp.~235--255.

\bibitem{HLFourier}
T.~P.~Hyt\"onen and M.~T.~Lacey,
\emph{Pointwise convergence of vector-valued Fourier series},
Math. Ann. \textbf{357} (2013), no.~4, 1329--1361.

\bibitem{HLWalsh}
T.~P.~Hyt\"onen and M.~T.~Lacey,
\emph{Pointwise convergence of Walsh--Fourier series of vector-valued functions},
Math. Res. Lett. \textbf{25} (2018), no.~2, 561--580.

\bibitem{HLP}
T.~P.~Hyt\"onen, M.~T.~Lacey, and I.~Parissis,
\emph{A variation norm Carleson theorem for vector-valued Walsh--Fourier series},
Rev. Mat. Iberoam. \textbf{30} (2014), no.~3, 979--1014.




\bibitem{HNVW}
T.~P.~Hyt\"onen, J.~M.~A.~M.~van Neerven, M.~C.~Veraar, and L.~Weis,
\emph{Analysis in Banach Spaces. Vol.~I. Martingales and Littlewood-Paley Theory},
Ergebnisse der Mathematik und ihrer Grenzgebiete, 3.~Folge, vol.~63,
Springer, Cham, 2016.

\bibitem{LaceyThiele}
M.~T.~Lacey and C.~Thiele,
\emph{A proof of boundedness of the Carleson operator},
Math. Res. Lett. \textbf{7} (2000), no.~4, 361--370.

\bibitem{Lepingle}
D.~L\'epingle,
\emph{La variation d'ordre $p$ des semi-martingales},
Z. Wahrscheinlichkeitstheorie Verw. Gebiete \textbf{36} (1976), no.~4, 295--316.

\bibitem{MSZ}
M.~Mirek, E.~M.~Stein, and P.~Zorin-Kranich,
\emph{Jump inequalities via real interpolation},
Math. Ann. \textbf{376} (2020), nos.~1--2, 797--819.

\bibitem{OSTTW}
R.~Oberlin, A.~Seeger, T.~Tao, C.~Thiele, and J.~Wright,
\emph{A variation norm Carleson theorem},
J. Eur. Math. Soc. \textbf{14} (2012), no.~2, 421--464.

\bibitem{Pisier}
G.~Pisier,
\emph{Martingales in Banach spaces (in connection with type and cotype)},
IHP lecture notes, February~9, 2011.

\bibitem{Pisier75}
G.~Pisier,
\emph{Martingales with values in uniformly convex spaces},
Israel J. Math. \textbf{20} (1975), nos.~3--4, 326--350.

\bibitem{PisierXu}
G.~Pisier and Q.~Xu,
\emph{The strong $p$-variation of martingales and orthogonal series},
Probab. Theory Related Fields \textbf{77} (1988), no.~4, 497--514.

\bibitem{Rubio}
J.~L.~Rubio de Francia,
\emph{Martingale and integral transforms of Banach space valued functions},
Probability and Banach Spaces (Zaragoza, 1985), Lecture Notes in Math., vol.~1221,
Springer, Berlin, 1986, pp.~195--222.

\bibitem{Sjolin}
P.~Sj\"olin,
\emph{An inequality of Paley and convergence a.e. of Walsh--Fourier series},
Ark. Mat. \textbf{7} (1969), no.~6, 551--570.

\bibitem{ThieleBook}
C.~Thiele,
\emph{Wave Packet Analysis},
CBMS Regional Conference Series in Mathematics, vol.~105,
American Mathematical Society, Providence, RI, 2006.

\bibitem{ThieleLectures}
C.~Thiele,
\emph{Lecture 5: Walsh phase plane; Lecture 6: Partial Walsh--Fourier sums},
lecture notes, 8 and 11 pp., respectively.

\bibitem{Weisz}
F.~Weisz,
\emph{Almost everywhere convergence of Banach space-valued Vilenkin--Fourier series},
Acta Math. Hungar. \textbf{116} (2007), nos.~1--2, 47--59.

\bibitem{Wenzel}
J.~Wenzel,
\emph{Mean convergence of vector-valued Walsh series},
Math. Nachr. \textbf{162} (1993), no.~1, 117--124.

\end{thebibliography}
\end{document}